\documentclass[authoryear,preprint,11pt,3p]{elsarticle}
\journal{Energy Policy}

\usepackage{graphicx}
\usepackage{xcolor}
\usepackage{physics}
\usepackage{amsmath}
\usepackage{amssymb}
\usepackage{amsthm}
\usepackage{eurosym}
\usepackage{makecell}
\usepackage{tabularx}
\usepackage{mhchem}
\usepackage{lscape}
\theoremstyle{plain}%

\theoremstyle{definition}

\usepackage{ulem}
\usepackage{nomencl}
\makenomenclature
\usepackage{longtable}
\makeatletter
\def\els@aparagraph[#1]#2{\elsparagraph[#1]{#2\@addpunct{.}}}
\def\els@bparagraph#1{\elsparagraph*{#1\@addpunct{.}}}
\makeatother
\usepackage{bm}
\usepackage{amsfonts}
\usepackage[colorlinks=True]{hyperref}
\usepackage{booktabs}
\usepackage{ragged2e}
\usepackage{tabularx}
\usepackage{multirow}
\usepackage{mathtools}
\usepackage[bb=boondox]{mathalfa}
\usepackage{cleveref}
\crefformat{section}{§#2#1#3}
\usepackage{nicefrac}
\usepackage{enumitem}
\usepackage{chngcntr}
\usepackage{booktabs}
\usepackage{tikz,tikz-3dplot}
\usepackage{pgfplots,pgfplotstable}
\usetikzlibrary{arrows,automata,calc,shapes,plotmarks, positioning,shapes.geometric,decorations.pathmorphing,patterns,fillbetween,pgfplots.groupplots}
\usepackage{circuitikz}
\pgfplotsset{compat=1.15}
\pgfplotsset{
    table/search path={plotdata},
}
\usepgfplotslibrary{colorbrewer,units}
\pgfplotsset{cycle list/Set3}
\definecolor{maincolor}{HTML}{032F99} %blue
\definecolor{secondcolor}{HTML}{ff5722} % dark orange
\definecolor{thirdcolor}{HTML}{c7d3d7}  %bluish gray
\definecolor{newtextcolor}{HTML}{bf360c}

\usepackage{adjustbox}

\usepackage{subcaption}

\usepackage{mathtools}

\usepackage{esvect}

\makeatletter
\newcommand{\oset}[3][0ex]{%
  \mathrel{\mathop{#3}\limits^{
    \vbox to#1{\kern-2\ex@
    \hbox{$\scriptstyle#2$}\vss}}}}
\makeatother

\begin{document}

\begin{frontmatter}
\title{Towards European Hydrogen Market Design: Perspectives from Transmission System Operators}

% Alternative title 
% European hydrogen market design: Perspective of Transmission System Operators

\author[1,2]{Marco Saretta\corref{cor1}}\ead{mcsr@dtu.dk}
\author[2]{Enrica Raheli}\ead{enri@ramboll.com}
\author[1]{Jalal Kazempour}\ead{jalal@dtu.dk}

\affiliation[1]{organization={Technical University of Denmark},
                city={Kgs. Lyngby},
                country={Denmark}}
\affiliation[2]{organization={Ramboll Denmark A/S}, 
                city={Copenhagen}, 
                country={Denmark}}
\cortext[cor1]{Corresponding author.}

\begin{abstract}
% Introduction to the topic – Context & Motivation
Despite hydrogen being central to Europe's decarbonisation strategy, only a small share of renewable hydrogen projects reached final investment decision. A key barrier is uncertainty about how future hydrogen markets will be designed and operated, particularly under Renewable Fuels of Non-Biological Origin requirements.
% Paragraph 2 – Research gap and method 
This study investigates the extent to which future hydrogen market design can be adapted from existing natural gas markets, and the challenges it must address. The analysis was based on a survey on future hydrogen market design targeting European gas transmission system operators, the stakeholder group responsible for operating the network and implementing the Network Codes. The survey provided an assessment of mechanism transferability and identified challenges for the early development of the hydrogen market.
% Preview findings - Part A: Transferability
Core market design principles and trading frameworks are broadly transferable from natural gas markets, as entry-exit systems and virtual trading points. Capacity allocation requires targeted adaptation to improve coupling with electricity markets. Tariffs require adaptation through intertemporal cost allocation, distributing infrastructure costs over time to protect early adopters. Balancing regimes should be revisited to reflect hydrogen's physical characteristics and different linepack flexibility usages.
% Preview findings - Part B: Challenges
Key challenges for early hydrogen markets include: temporal mismatches between variable renewable supply and expected relatively stable industrial demand, limited operational flexibility due to scarce storage and reduced pipeline linepack, fragmented regional hydrogen clusters, and regulatory uncertainty affecting long-term investment decisions.
% What your work teaches us
These findings provide empirical input to the hydrogen Network Codes and offer guidance to policymakers designing hydrogen market frameworks.

% 5. What your work teaches us
\end{abstract}
%%Graphical abstract
% \begin{graphicalabstract}
% \includegraphics{grabs}
% \end{graphicalabstract}


\end{frontmatter}
\section*{Nomenclature}
\vspace{-4mm}
{\small
\addtocounter{table}{-1}  % counteracts the count after the table
\begin{longtable}{p{3.1cm} p{12.5cm}}\label{tab:acronyms}\\
\multicolumn{2}{l}{\textbf{Acronyms}}  \\
ACER  & Agency for the Cooperation of Energy Regulators \\
DSO   & Distribution System Operator \\
EEX   & European Energy Exchange AG \\
EHB   & European Hydrogen Backbone \\
ENTSOE & European Network of Transmission System Operators for Electricity \\
ENTSOG & European Network of Transmission System Operators for Gas \\
ENNOH & European Network of Network Operators for Hydrogen \\
EU & European Union \\
\ce{H2} & Hydrogen \\
OTC   & Over-the-Counter \\
RAB   & Regulatory Asset Base \\
NC & Network Codes \\
NC BAL & Network Code on Balancing \\
NC CAM & Network Code on Capacity Allocation Mechanisms \\
NC TAR & Network Code on Tariff Structures \\
RFNBO & Renewable Fuels of Non-Biological Origin \\
STP & Standard Temperature and Pressure \\
TSO   & Transmission System Operator \\
TTF   & Title Transfer Facility \\
VTP   & Virtual Trading Point \\
WDO   & Within-Day Obligation \\
\end{longtable}}

\section{Introduction}
\label{sec:introduction}

\begin{table}[bp]
\centering
\caption{Key physical and operational differences between hydrogen and natural gas transport. Density reported at Standard Temperature and Pressure (0°C, 1 bar).}
\label{tab:ng_vs_h2_pipeline_short}
\begin{tabular}{@{}lll@{}}
\toprule
\textbf{Parameter} & \textbf{Natural Gas} & \textbf{Hydrogen} \\ \midrule
Density (STP) & $\sim$0.75 kg/\ce{m^3} & 0.09 kg/\ce{m^3} \\
Lower Heating Value & 50 MJ/kg & 120 MJ/kg \\
Volumetric Energy Density & 39 MJ/\ce{m^3} & 10.8 MJ/\ce{m^3} \\
Linepack Capacity & High & Low \\
Network Maturity & Pan-European & Fragmented (early phase) \\
Storage Availability & Mature & Limited (early phase) \\
\bottomrule
\end{tabular}
\end{table}

\subsection{Context and research gap}
% Hydrogen's critical role and EU ambitions
As Europe transitions towards climate neutrality, low-carbon and renewable hydrogen have been identified as essential for decarbonising hard-to-abate sectors \citep{Shen2024TheImpacts}. Hydrogen's central role is also reflected in ambitious European targets and support schemes: the REPowerEU Strategy aims to reach 10 million tonnes of renewable hydrogen production by 2030, supported by 40 GW of electrolysers,  and dedicated financial schemes \citep{EuropeanHydrogenObservatory2025EUDeal}. Furthermore, the European Hydrogen Backbone initiative aims to develop the pan-European infrastructure and pipelines needed to connect producers and consumers, led by 33 future hydrogen Transmission System Operators (TSOs) \citep{EuropeanHydrogenBackbone2024EuropeanCompetitiveness}. Finally, the European Commission has published a Delegated Act on Renewable Fuels of Non-Biological Origin (RFNBO) \citep{EuropeanParliament2023DelegatedOrigin}, which defines requirements for low-carbon and renewable hydrogen producers and mandates strict production constraints, including temporal and geographical correlation with renewable generation, additionality, and lifecycle emission thresholds. The strictness of these constraints carries a measurable cost.  \cite{Brandt2024CostFramework} show that allowing unrestricted grid-electricity use in hydrogen production could lower the on-site hydrogen supply cost by \euro  0.55-1.62/kg \ce{H2} relative to the compliant pathway, illustrating the trade-off regulators face between environmental integrity and commercial viability.

% Gap: Projects stall despite policy support
Meeting this ambition requires a very large-scale investment. \citep{Nunez-Jimenez2022CompetitiveUnion} estimate that supplying the projected EU hydrogen demand would require \$2-2.4 trillion in renewable generation, electrolysis, and transport infrastructure.  Yet ambition has not translated into investment: despite these targets, only 25\% of the expected 2030 renewable hydrogen capacity has reached final investment decision \citep{Fraile2025Clean2025}. One reason for this is the uncertainty around future market design. This uncertainty is not purely technical but also political. \cite{Guinlle2026CompetingMarket} discuss that EU hydrogen governance has left unresolved a persistent conflict between incumbent actors favouring blending into existing gas infrastructure (`system preservation') and new entrants requiring dedicated renewable hydrogen infrastructure (`system creation'), with neither pathway backed by a binding enforcement mechanism. Without knowing the rules under which the market will operate, developers and off-takers are negatively affected in the bankability of their projects \citep{Lagioia2023BlueAffairs}.

\begin{figure}
    \centering
    \includegraphics[width = 0.9\linewidth, trim={0cm 7.5cm 0cm 7cm},clip]{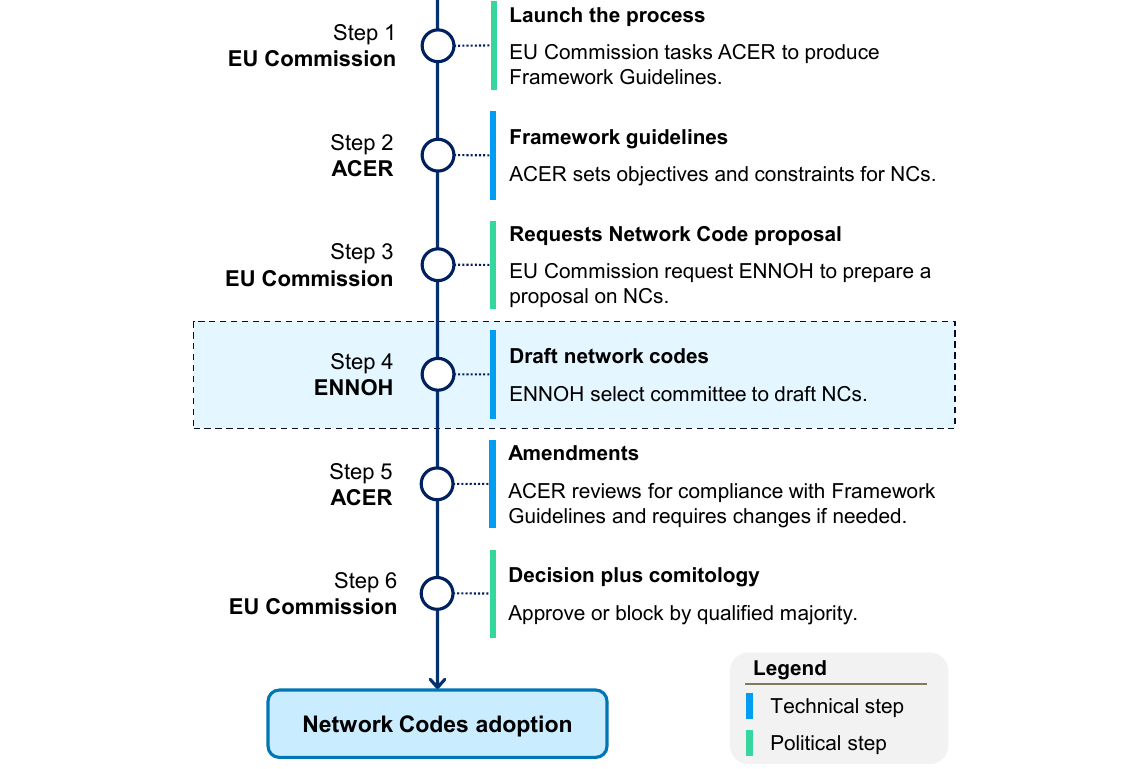}
    \caption{Network Codes development process for hydrogen in the European Union. This publication aims to contribute to the fourth step by providing informed insights useful to ENNOH in drafting the Network Codes. Image adapted from \citet{ENTSOE2023ProcessGuidelines}.}
    \label{fig:nc_developement}
\end{figure}

These rules are expected to be codified in dedicated Network Codes (NC). However, the development of NC is often a lengthy, multi-party and multi-stage regulatory process, illustrated in Figure \ref{fig:nc_developement}, in which the drafting responsibility is assigned to the European Network of Network Operators for Hydrogen (ENNOH) under the supervision of the Agency for the Cooperation of Energy Regulators (ACER). Until the hydrogen NC enters into force, investors face an uncertain market framework, discouraging investment in projects.

% Physical differences prevent gas market replication
While the EU natural gas market is often referenced as a model for future hydrogen market design \citep{Johnson2025RealisticTransition}, fundamental physical and operational differences (Table \ref{tab:ng_vs_h2_pipeline_short}) prevent direct replication. Hydrogen's density is approximately one-tenth that of natural gas, directly reducing the mass of gas that can be stored in pipelines at a given pressure. Although hydrogen has a higher lower heating value, its lower density results in a substantially lower volumetric energy density, reducing available linepack flexibility (i.e., the energy stored in pipelines for balancing purposes). From an operational perspective, infrastructure limitations aggravate these physical constraints. Early European networks are expected to emerge as fragmented regional clusters rather than interconnected systems, and will initially lack mature bulk storage infrastructure. Hence, future market design cannot simply replicate gas market mechanisms, and it must account for hydrogen distinct physical and operational characteristics.

% Literature review: design and regulatory perspectives
Recent research addressed hydrogen market design in Europe from multiple perspectives. A central reference is \citet{Steinbach2024TheTrading}, which derives three key market design criteria through interviews with stakeholders across the hydrogen value chain: policy support measures, infrastructure regulation, and trading arrangements. The study highlights the need in these areas to reduce regulatory uncertainty, provide targeted financial support, and introduce fundamental natural gas market principles, such as third-party access, unbundling, and entry-exit systems. However, the analysis remains at a conceptual level and does not specify how these elements translate into market design. Structured design analyses complement this perspective, with the work from \cite{Niedrig2024MarketMarket} examining how market design evolves across development stages. The study finds that, in the early phases, a fragmented infrastructure and limited participation are expected to favour bilateral contracts, supporting investment but constraining liquidity. As networks expand and market areas integrate, exchange trading becomes feasible, improving transparency and short-term efficiency. From a regulatory economics perspective, \citet{Martinez-Rodriguez2026PipelineNetworks} argue that hydrogen markets differ from natural gas markets due to the need for new infrastructure and weaker natural-monopoly characteristics. The study proposes a model based on negotiated access, pipeline-level regulation, and long-term capacity contracts to support investment, complemented by proportionate regulation and market-based coordination mechanisms. Finally, broader policy syntheses situate these design questions in a cross-sectoral context. \citet{Hasan2025DesigningReform} apply a Structure-Conduct-Performance-Regulation framework, informed by consultations with over 40 stakeholders. This derives recommendations on pricing mechanisms, business models, and infrastructure access, drawing on earlier developments in the electricity, natural gas, and liquefied natural gas (LNG) markets. The study identifies market reference prices, government-backed offtake agreements, and third-party access provisions as the principal enablers of market formation. However, as with \citet{Steinbach2024TheTrading}, the recommendations remain at the level of policy instruments and are not translated into the operational rules governing network access and system operation.

% Literature review: physical boundary conditions
Physical and system-level studies define the boundary conditions within which these designs must operate. \citet{Sargent2025LinepackModelling} show that hydrogen provides significantly less usable linepack than natural gas, limiting the feasibility of balancing regimes based on passive flexibility. Furthermore, integrated energy system models confirm the long-term role of hydrogen networks in Europe \citep{Neumann2023TheEurope}, particularly in linking renewable supply regions with demand centres and reducing system costs, but typically assume perfect coordination and do not account for detailed market rules, such as access and balancing mechanisms.

% Research gap
The literature converges on several conclusions: regulatory clarity is a prerequisite for investment, gas market mechanisms are only \textit{partially} transferable to hydrogen, and physical constraints narrow the operational flexibility available to market designers. Despite these research efforts, a gap remains at the granularity at which Network Codes are drafted: existing contributions either stop at design principles and policy instruments, or examine a single regulatory area in isolation. The operators who will run these networks, and whose technical input feeds into ENNOH's drafting mandate, are routinely consulted in regulatory processes, but no published study has assembled their positions as a distinct group across countries and these areas jointly. This study addresses that gap.

\subsection{Research questions and contributions}
This paper addresses the identified gap via two research questions: (\textit{i}) to what extent can future hydrogen design be adapted or replicated from the current EU natural gas one? and (\textit{ii}) what challenges must the hydrogen market design account for? To address these questions, a structured survey was developed, targeting European gas TSOs. Hydrogen market design will affect producers, industrial consumers, storage operators, traders, regulators, and electricity market operators. This study consults one of these families, the TSOs, as the group from which ENNOH is constituted and which therefore supplies the technical input to the drafting of the Network Codes. The resulting view is therefore partial. The survey results lead to two main contributions, both aimed at informing ENNOH during the drafting of hydrogen Network Codes in Step 4 of Figure \ref{fig:nc_developement}. First, it assesses the transferability of gas market mechanisms to hydrogen across the five regulatory areas, identifying which are replicable and which require minor or major adaptations due to technical and operational constraints. Second, it presents a risk analysis of challenges from TSO consultation, organised into operational and regulatory challenges.

\subsection{Paper structure}
The paper is structured as follows. Section \ref{sec:methods} presents the survey design and data analysis. Sections \ref{sec:transferability} and  \ref{sec:challenges} present the results of the study, respectively, assessing the transferability of gas market components to hydrogen and identifying the challenges in designing future hydrogen markets. Section \ref{sec:discussion} discusses implications for market design and infrastructure policy. It is assumed that the reader is familiar with the general concepts of EU gas markets, e.g., entry-exit systems, capacity allocation, nominations, balancing frameworks, virtual trading points, and commodity trading. For readers unfamiliar with these concepts, a comprehensive technical overview is provided in the supplementary materials \citep{Saretta2026DesigningRepository}.
\section{Methods}
\label{sec:methods}

\begin{figure}
    \centering
    \includegraphics[width = \linewidth, trim={0cm 11.75cm 0cm 10cm},clip]{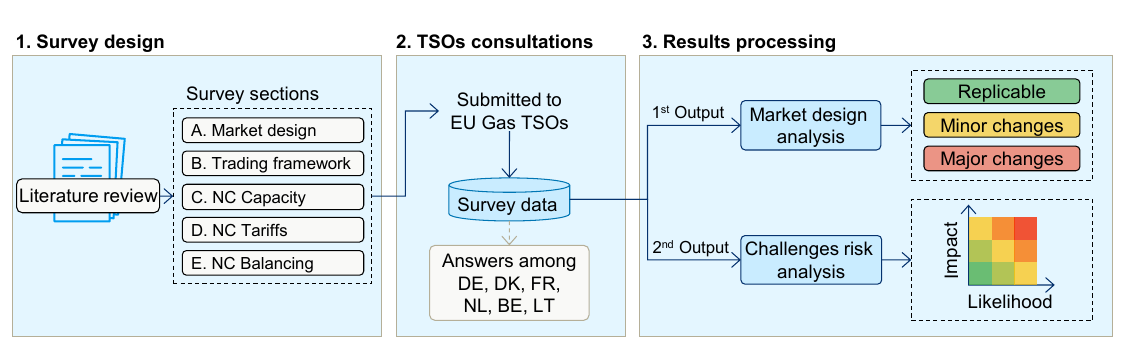}
    \caption{Three-stage methodology for assessing hydrogen market design through TSO consultation.}
    \label{fig:methodology}
\end{figure}

To address research questions (i) and (ii), the study employed a survey to collect the insights and perspectives from EU gas TSOs. Within this study, TSOs were selected as the survey recipient  among the broader set of hydrogen market stakeholders, which includes producers, industrial consumers, storage operators, traders, regulators, and electricity market operators. This choice reflects their operational expertise in gas network management, their roles in hydrogen infrastructure development, and their expected participation in ENNOH, the body responsible for drafting the hydrogen Network Codes. The underlying methodology was organised into three stages (Figure \ref{fig:methodology}): survey design, TSO consultation, and results processing. Complementary details on the methodology are provided in \ref{sec:survey}. 

\subsection{Survey design}

A survey was developed to assess the transferability of EU gas market mechanisms to hydrogen networks and identify future challenges, totalling 37 questions. The survey was organised into five sections: (a) Market and Regulatory Design, (b) Trading Framework, (c) Network Code on Capacity Allocation Mechanisms (NC CAM), (d) Network Code on Tariff Structures (NC TAR), and (e) Network Code on Balancing (NC BAL). Sections (a) and (b) addressed general market design and trading mechanisms derived from the EU Hydrogen and Decarbonised Gas Market Package \citep{EuropeanCommission2025HydrogenMarket}, while sections (c), (d), and (e) explored the adaptation of current EU natural gas network codes for hydrogen. Each section had its questions organised into subtopics. Section (a) covered entry-exit systems, hydrogen day/year, the role of shippers and Distribution System Operators (DSO) cooperation. Section (b) included questions on Virtual Trading Points (VTPs), Over-The-Counter (OTC) trading, the role of exchange-based platforms and general pricing structure. Section (c) addressed capacity products, capacity allocation mechanisms, nomination and renomination procedures, building from the natural gas NC CAM \citep{Entsog2017CapacityCode}. Section (d) covered capacity pricing and congestion procedures based on the natural gas NC TAR \citep{Entsog2018TariffOverview}. Section (e) addressed balancing regimes, the role of linepack flexibility, and within-day obligations (WDOs), based on the natural gas NC BAL \citep{Entsog2018BalancingOverview}. Questions across the five sections were formulated as a mix of open-ended and yes/no formats. Before its distribution, the survey was pilot-tested with Energinet (the Danish TSO) in August 2024, after which the questions were refined based on the feedback received. 

\subsection{TSO consultation and data collection}
All European gas TSOs are listed in the ENTSOG member directory, distinguishing three categories: members, associated partners, and observers. The survey was addressed to the 43 members, as they are the certified TSOs to which the EU Network Codes and the Hydrogen and Decarbonised Gas Market Package apply, and from which ENNOH is constituted. Associated partners operate outside this legal framework, while observers are not certified TSOs, so neither group is governed by the arrangements whose transferability this study assesses. Each member was approached through a named contact person, otherwise, the survey was sent through the contact form or information address published on the operator’s website, or through the operator’s general corporate address. Members were contacted between September 2024 and February 2025.

To preserve confidentiality, the TSO organisations and individual respondents were anonymised, with responses reported only at the country level. Excluding non-responses and declines to participate, the survey received 7 responses from TSOs across 6 countries: Belgium, Denmark, France, Germany, Lithuania, and the Netherlands. These countries account for approximately 45\% of the planned length of the European hydrogen pipeline by 2040, according to the European Hydrogen Backbone data \citep{ENTSOG2026HydrogenMap}. 
% R2C1
They also account for 50\% of EU-27 hydrogen production and consumption in 2024 across all production routes \cite{EuropeanHydrogenObservatory2026HydrogenProduction}, for 67\% of binding hydrogen offtake commitments in the European Economic Area as of May 2025 \cite{ACER2025EuropeanReport}, and for 57\% of European electrolyser capacity in operation and under construction \cite{EuropeanHydrogenObservatory2025TheLandscape}. Table \ref{tab:respondent_profile} characterises each participating country in terms of network maturity, regulatory status, and planned infrastructure. Interviewees held positions in operations management, strategic planning, and regulatory affairs.

\begin{table}[h]
\centering
\footnotesize
\caption{Profile of the participating countries at the time of the survey (September 2024 to February 2025).}
\label{tab:respondent_profile}
\begin{tabularx}{\textwidth}{@{}>{\bfseries\RaggedRight}X>{\RaggedRight}X>{\RaggedRight}X>{\RaggedRight}X>{\RaggedRight}X@{}}
\toprule
Country & Region & H$_2$ network maturity & Regulatory status & Planned infrastructure \\
\midrule
Germany & Central & Core network approved, commissioning to 2032 & Access, balancing and amortisation account consultations published & National core network, about 60\% repurposed \\
\addlinespace
Netherlands & North-West & National network under construction & Interim tariff set by ministry; conditions published & Phased national rollout \\
\addlinespace
Belgium & North-West & First segments under construction & Network operator designated; first ten-year plan published & Backbone linking clusters and import terminals \\
\addlinespace
France & Central & Cluster-based; no national backbone & No published tariff or access rules at the time & Cross-border interconnection; H2Med \\
\addlinespace
Denmark & Northern & Backbone investment decision taken & Consultations begun on planning and tariff principles & First backbone segment \\
\addlinespace
Lithuania & Baltic & No network; feasibility stage & Strategic guidelines only; transposition pending & Nordic-Baltic Hydrogen Corridor \\
\bottomrule
\end{tabularx}
\end{table}

\subsection{Data processing and synthesis}
Survey data was processed as follows. Responses were coded by component. Answers indicating no position, such as ``no information'' or ``under discussion'', were recorded as unresolved and treated as neither agreement nor dissent. Each remaining answer was classified according to the change it implied relative to the corresponding gas arrangement: none, a change to a parameter or procedure within it, or a structural change. Each component was assigned the most demanding level recorded. Two further cases were also classified as requiring major modifications: components on which respondents described different national designs, and components for which no settled arrangement exists in any national framework, since in both cases the gas arrangement cannot be carried over without a design decision that has not yet been taken. No majority threshold was applied. Response coverage was uneven: between three and six of the seven respondents recorded a position on each component, with the trading framework the thinnest. The number of substantive responses is stated in the text wherever fewer than five respondents addressed a component. Statements attributed to surveyed TSOs report survey responses, where statements introduced as our assessment are the authors' interpretation and are not attributed to respondents.

\section{Transferability of gas market mechanisms to hydrogen}
\label{sec:transferability}

\begin{figure}
\centering
\includegraphics[width = \linewidth, trim={0cm 11.cm 0cm 9.8cm},clip]{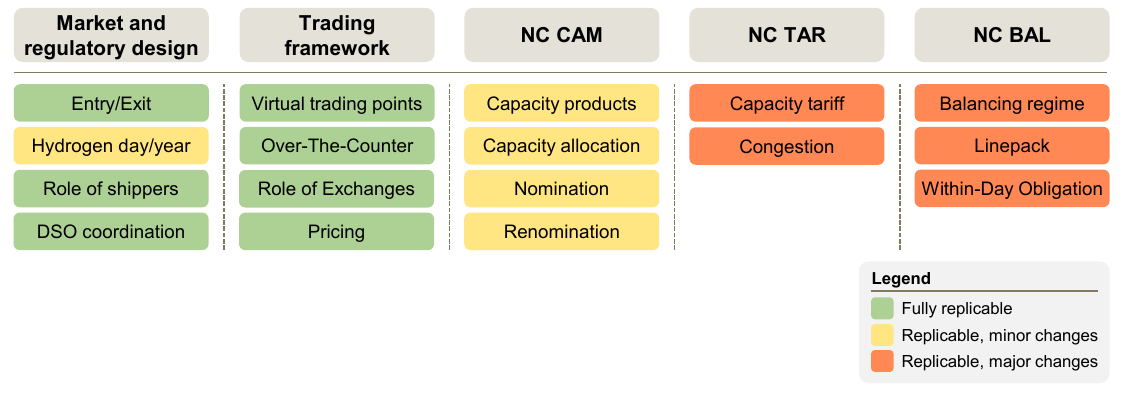}
\caption{Transferability assessment of natural gas market design for the hydrogen market. Market design components across five regulatory areas are classified into three categories following the rule in Section 2.3: fully replicable (green), replicable with minor changes (yellow), and replicable with major changes (orange). Classification is at the component level; where greater nuance is required, this is provided in the corresponding section. The assessment based on consultations with seven European gas TSOs from Belgium, Denmark, France, Germany, Lithuania, and the Netherlands.}
\label{fig:h2_market_properties}
\end{figure}

This section assesses the extent to which the current gas market can serve as a reference model for the design of the future hydrogen market. The assessment has been conducted for each of the five survey sections, namely market and regulatory design, trading framework, Network Codes on capacity allocation mechanisms, Network Codes on capacity tariffs, and Network Codes on balancing schemes. Each survey section subtopic has been mapped to a marked design component and further categorised as fully replicable, replicable with minor modifications, or replicable with major changes, as displayed in Figure \ref{fig:h2_market_properties}. 

The five areas are not at the same regulatory stage: the entry-exit system and the network user entity are mandated by EU Hydrogen and Decarbonized Gas Market Package \cite{EuropeanCommission2024EUPackage} and are reported here as legislated outcomes, whereas the hydrogen day and year, capacity products, allocation, gate closures, linepack and the balancing regime are not fixed at EU level, and it is there that the divergence reported below arises.
\subsection{Market and regulatory design}
\label{subsec:market_regulatory}

The market and regulatory design survey section was structured around four subtopics: entry-exit systems, hydrogen day and year, shipper roles and evolution, and DSO coordination. Surveyed TSOs indicated that the transferability of these components is largely determined by the EU Hydrogen and Decarbonised Gas Market Package (Directive 2024/1788 and Regulation 2024/1789) \citep{EuropeanCommission2025HydrogenMarket}.

% entry / exit and h2 day and year
The entry/exit system is mandated for hydrogen networks by the EU Hydrogen and Decarbonised Gas Market Package from January 2033. Under this system, a network user would book capacity at entry and exit points independently, with hydrogen injected or withdrawn respectively at entry and exit points in the network. The Package, however, does not define the hydrogen day and year, as these need to be decided by the member states. Five of the six participating TSOs that addressed this point indicated that the hydrogen day and year should be aligned with the calendar day (00:00-24:00 daily, January-December annually), moving away from the gas market conventions (06:00-06:00 UTC gas day, October-September gas year). The remaining respondent, in a market at feasibility stage, expected the gas structure to be carried over. TSOs emphasised that this shift from gas day/year to calendar day/year represents a necessary adjustment. The traditional gas market day and year structures were originally designed to match daily and seasonal heating demand patterns. In contrast, aligning the hydrogen day with the calendar day facilitates coupling with electricity markets, as hydrogen demand is expected to follow industrial load profiles and its production is tightly coupled to the electricity market.

\begin{figure}
    \centering
    \includegraphics[width = \linewidth, trim={0cm 12cm 0cm 10.5cm},clip]{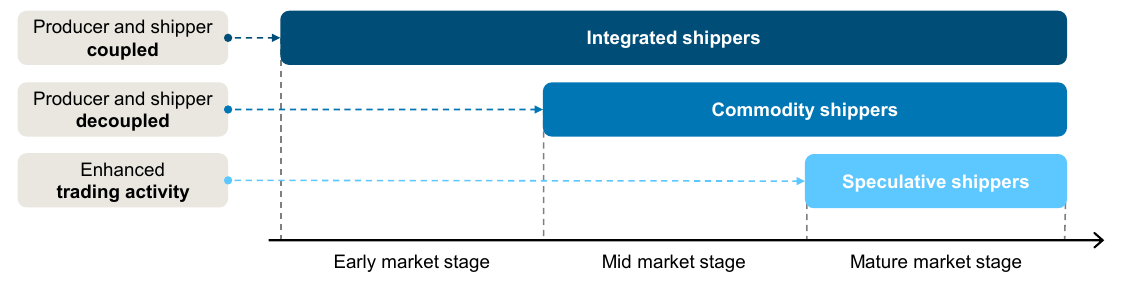}
    \caption{Expected evolution of shipper types across hydrogen market development stages according to surveyed TSOs.}
    \label{fig:shipper_evolution}
\end{figure}

% Shipper role
The Hydrogen and Decarbonised Gas Market Package also indicates maintaining the entity of network users, also known as shippers, inherited from existing natural gas market structures. Network users are defined as customers or potential customers of a system operator involved in transport and balancing functions, who are responsible for balancing their hydrogen injections and withdrawals \citep{EuropeanCommission2025HydrogenMarket}. The terms shipper and network user are used interchangeably throughout this paper. Surveyed TSOs expect the shipper's role to evolve along different market phases. As the hydrogen market evolves from early low-liquidity stages to more mature phases, different shipper types are expected to emerge, as illustrated in Figure \ref{fig:shipper_evolution}. In the early stages, integrated shippers are expected to dominate: these are large energy companies that handle both production and consumption, managing the entire supply chain from production to delivery with limited trading interactions. As hydrogen availability and demand grow, commodity shippers are expected to emerge, focusing on buying hydrogen from producers and selling it to consumers without direct involvement in production. Finally, in mature market phases characterised by high liquidity, speculative shippers and traders are anticipated to engage primarily in arbitrage, leveraging price differentials across virtual trading points on energy exchanges. Speculative shippers are the last to emerge, as their activity requires a high volume of trades and sufficient market liquidity to generate profitable arbitrage opportunities. Two respondents described the typology in detail, and a third endorsed it rather than giving an independent account. 

% DSO coordination
TSOs anticipated that DSO cooperation frameworks would be maintained, though DSO involvement would depend on the development of the low-pressure distribution infrastructure. Hydrogen networks are expected to be deployed first through high-pressure transmission systems and later through low-pressure distribution, as liquidity grows. Unlike natural gas DSOs, which mainly serve households and smaller industries, hydrogen is expected to serve predominantly large industrial consumers. Therefore, DSO involvement will likely be limited in the early stages, as it is expected that TSOs will handle hydrogen transmission until sufficient liquidity and infrastructure develop.
\subsection{Trading framework}
\label{subsec:trading}

The trading framework section in the survey was structured around four subtopics: virtual trading points, Over-The-Counter trading, exchange-based platforms, and pricing mechanisms. Response coverage was thin across this section, and the respondents who did address these components, between three and four of the seven depending on the component, expected them to be transferable from gas markets.

Three of the seven respondents addressed virtual trading points (VTPs), and described them as foundational for a functional hydrogen market: without them, hydrogen can only be traded through physical delivery. VTPs enable the title transfer of gas ownership, while hydrogen remains within pipelines, allowing network users to trade bilaterally and independently of physical location within the entry-exit zone.

Four respondents addressed whether commodity trade would follow Over-The-Counter (OTC) trading, exchange-based platforms, or both, and all four expected OTC trading first, with exchange trading following market maturity. OTC trading involves bilateral contracts in which the parties agree on the volume of hydrogen to be exchanged, its price, and the contract's tenor or duration. OTC agreements are non-standard contracts because they are customised through negotiations between the parties, and are typically non-public. Exchange-based trading occurs on organised platforms with standardised contracts and public price discovery, requiring VTPs as underlying infrastructure. The hydrogen market is anticipated to begin as OTC-dominated due to initial low liquidity, with producers and importers securing long-term supply contracts directly with consumers via bilateral agreements. As networks interconnect and liquidity increases, exchange-based trading is expected to develop alongside OTC markets, using VTPs for continuous trading.

Three respondents addressed pricing mechanisms, and expected them to follow exchange-based prices as the market matures. In the early stages, it is expected that OTC trades will dominate and, consequently, prices will be less transparent, as OTC contract details are not public. As liquidity increases, continuous price signals are expected to emerge at trading exchanges facilitated by VTPs, analogously to how the Dutch Title Transfer Facility (TTF) became the reference VTP for natural gas trading in the Netherlands. Exchange trading would likely represent only a portion of total volumes, with the remainder traded OTC, though exchange prices would serve as market reference points.
\subsection{Network Code on Capacity Allocation Mechanisms}
% intro
The survey section on capacity allocation included four subtopics: capacity products, capacity allocation mechanisms, nomination and renomination procedures. Surveyed TSOs broadly indicated that these elements are transferable from gas markets with minor adjustments.

% Capacity products
For capacity products, TSOs indicated a preference for issuing them as long-term contracts in early market stages, typically as yearly or multi-year agreements, because substantial infrastructure investment requires network users to demonstrate long-term commitment. Contract lengths vary among surveyed TSOs. The Danish TSO reported plans to offer capacity products ranging from short-term to one-year annual contracts for up to 15 years ahead, with a target booking level of 0.5 GW for at least 10 years needed to justify infrastructure investment \citep{Energinet2025CapacityBackbone}. In Germany, the Bundesnetzagentur (federal network regulator) established frameworks that propose standard yearly capacity products while reserving at least 10\% of capacity for shorter-term products, such as monthly and daily contracts \citep{Bundesnetzagentur2024SecondWaKandA}. A minor change from gas markets would be to align yearly products with calendar years rather than gas years, reflecting hydrogen's industrial demand patterns rather than seasonal heating cycles.

% capacity allocation
Regarding capacity allocation methods, surveyed TSOs generally indicated a preference for the first-come first-served allocation method in early market phases, representing a divergence from natural gas markets where auction mechanisms predominate. This preference is due to the expected limited number of market participants during initial network development and consequent low liquidity. German frameworks specify that when demand exceeds available capacity, allocation should transition to auction-based mechanisms to ensure efficient distribution \citep{Bundesnetzagentur2024SecondWaKandA}.

Surveyed TSOs agreed that nomination and renomination procedures should be maintained, although gate-closure structures require adaptation to enable coupling with the electricity market. Four of the seven respondents addressed nomination gate closure specifically; at the time of the survey, only the Danish and German TSOs had developed detailed proposals, with the remaining respondents indicating that these aspects were still under development.
% First Energinet approach
The Danish TSO proposal is described as follows. It proposes opening hydrogen nominations after the day-ahead electricity market has cleared, allowing shippers to submit nominations based on the realised market outcomes \citep{Energinet2024ProposalNetwork}. This is particularly relevant for market participants whose production is closely coupled to electricity markets, e.g., electrolyser-based producers. To incentivise accurate nominations, it is also proposed to implement a fee for nomination errors, as accurate nominations would reduce unforeseen hydrogen flows and system imbalances, enabling the TSO to minimise reserve linepack requirements and maximise available linepack flexibility for network users. The specific nomination time unit for Denmark, whether 15-minute intervalsaligned with electricity markets or 1-hour blocks, had not been finalised at the time of the survey.
% Germany approach
The German proposal is contained in the WaKandA document \citep{Bundesnetzagentur2024SecondWaKandA}. It states that shippers are required to nominate hourly entry and exit quantities, with TSOs defining uniform deadlines aligned with balancing requirements. Exit nominations are required only for storage injection, cross-border transport, or where multiple shippers book the same exit point across different balancing groups. A key difference between the German and Danish proposed nomination procedures lies in the specificity of the timing. Denmark anchors gate closures to day-ahead market outcomes as a design principle, while Germany delegates specific times to TSO coordination, with implementation expected by October 2026.

% Renomination
Renomination procedures would allow shippers to adjust nominations in response to intraday electricity price signals. The Danish TSO indicated that network users must be able to make renominations, which are allowed immediately after the first nomination deadline, while more specific details are still under development \citep{Energinet2024ProposalNetwork}. German frameworks allow multiple intraday renominations without frequency limits, provided that the deadlines allow sufficient time for TSO processing \citep{Bundesnetzagentur2024SecondWaKandA}.

\begin{figure}
    \centering
    \includegraphics[width = \linewidth, trim={0.5cm 9cm 0.5cm 12cm},clip]{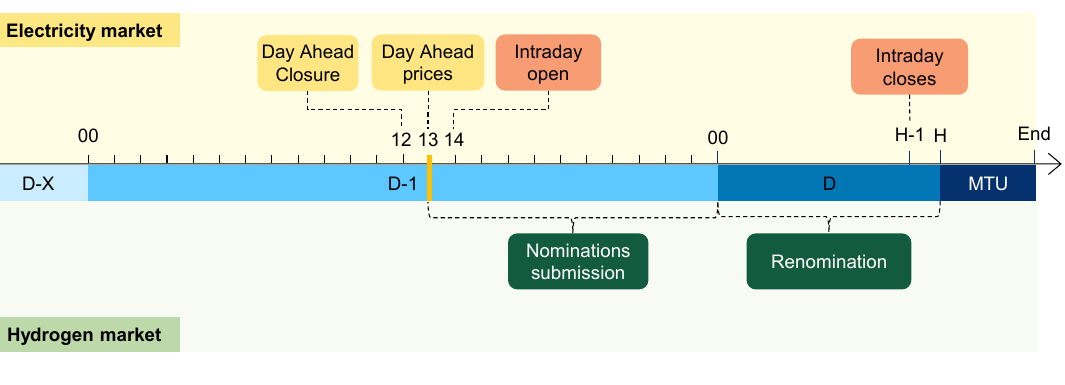}
    \caption{Illustrative alignment between electricity market timing and hydrogen nomination. The electricity side reflects current market coupling, which has operated on a 15-minute market time unit in the day-ahead timeframe since October 2025. The hydrogen side is a schematic based on the timings provided by surveyed TSOs. The Figure aims to convey the intuition behind aligning nomination with day-ahead clearing and renomination with intraday trading, and is not intended as a definitive gate-closure schedule.}
    \label{fig:gate_closure}
\end{figure}

Figure \ref{fig:gate_closure} sets hydrogen nomination against the electricity timetable. While the electricity side is fixed, with day-ahead coupling cleared on a 15-minute resolution since October 2025, the hydrogen side is not. No regulation currently governs hydrogen nomination gate closures, but surveyed TSOs expect them to align with the electricity market, and Figure \ref{fig:gate_closure} illustrates the intuition behind that alignment rather than a proposed schedule. The motivation follows from the production schedule of electricity-coupled hydrogen producers: an electrolyser's output decision follows the day-ahead price, so opening nomination after day-ahead clearing lets a shipper nominate against a realised price rather than a forecast, while a renomination window after intraday trading lets the hydrogen schedule track the electricity schedule. What remains open is the granularity, with the German respondent describing metering values every 15 minutes within an hourly balancing period, and the French respondent expecting nomination timelines to shorten to between 15 minutes and one hour. Furthermore, the coupling reaches beyond nomination timing. Renewable variability passes through the electrolyser as a within-day imbalance, with linepack buffering far less than in gas systems, and surveyed TSOs reported balancing periods shortened to below one hour in response. They did not expect binding congestion in the early years, given the limited number of participants. In our assessment, the flexibility an electrolyser can offer the power system depends on how much imbalance the balancing regime allows it to carry, coupling the two markets through linepack allocation rather than price alone.
\subsection{Network Code on Tariff Structures}
\label{subsec:nc_tar}

% i have to remember to mention in the graph that 
Network tariff design in the survey covered two components: capacity tariff
structures and congestion management. Surveyed TSOs indicated that both would
require fundamental redesign to avoid penalising early market players.

\begin{figure}
\centering
\begin{subfigure}{0.48\textwidth}
    \includegraphics[width=\linewidth, trim={0cm 11cm 9.5cm 9cm},
    clip]{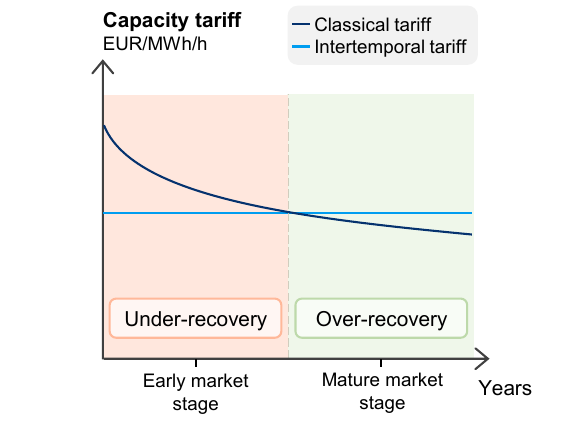}
    \caption{Classical and intertemporal tariff schemes comparison.}
    \label{fig:inter_1}
\end{subfigure}
\hfill
\begin{subfigure}{0.48\textwidth}
    \includegraphics[width=\linewidth, trim={9.5cm 11cm 0cm 9cm},
    clip]{Figures/capacity_tariffs.pdf}
    \caption{Projected regulatory debt under intertemporal tariff.}
    \label{fig:inter_2}
\end{subfigure}
\caption{Intertemporal cost allocation under hydrogen entry-exit regulation.
Panel (a) shows classical versus intertemporal tariffs across market
stages. Panel (b) shows regulatory debt evolution over time. Adapted from
\cite{Martinez-Rodriguez2026PipelineNetworks}.}
\label{fig:intertemporal_cost}
\end{figure}

% Survey outcome to be streamlined
Capacity tariffs, generally expressed in EUR/MWh/h, grant network users the right to inject or withdraw hydrogen up to a contracted capacity level and constitute a primary revenue stream for TSOs. These revenues are regulated to ensure cost recovery of the Regulatory Asset Base (RAB), which reflects the value of approved network assets, including hydrogen pipelines and repurposed gas infrastructure. For hydrogen markets, a central challenge arises: if the classical tariff scheme from natural gas markets were applied, RAB-related costs would need to be recovered across low transported volumes, resulting in high per-unit tariffs. Such high tariffs would be borne by early network users, potentially deterring participation and slowing demand growth. Therefore, directly adopting traditional natural gas tariff schemes for hydrogen could threaten market development at a stage when network uptake is most critical. 

% Intertemporal cost allocation
To address this, the intertemporal cost allocation scheme has been introduced in Article 5 of EU Regulation 2024/1789 \citep{EuropeanParliament2024RegulationHydrogen}, which allows TSOs to set tariffs below full cost-recovery levels in early market stages, thereby postponing the resulting revenue gaps to later periods, once network utilisation has grown. Figure \ref{fig:intertemporal_cost} illustrates the two core elements of the mechanism. Panel \ref{fig:inter_1} displays illustrative classical and intertemporal tariff projections. Under the classical approach, tariffs are initially high because RAB costs are recovered across low initial transported volumes, and they decline as market liquidity grows and more capacity is booked. In contrast, intertemporal tariffs are set at a lower level in early stages, resulting in initial under-recovery, followed by over-recovery as market liquidity rises, with the classical tariff eventually falling below it. While shown as constant for illustrative purposes, in practice, the intertemporal tariff is periodically adjusted as demand forecasts and cost projections are updated. Panel \ref{fig:inter_2} shows the resulting regulatory debt. During the under-recovery phase, uncollected revenues accumulate as debt, which peaks at the transition to over-recovery. As revenues then exceed allowed costs, the debt is gradually repaid until the break-even point is reached. The mechanism, therefore, shifts cost recovery forward in time while ensuring that network users ultimately bear all infrastructure costs. For a comprehensive description of the mechanism, the authors refer to \citet{Martinez-Rodriguez2026PipelineNetworks}. Six respondents addressed capacity tariffs: one described the intertemporal cost allocation mechanism directly, as the approach that its national framework would apply; one noted that it was under discussion at the European level; the remaining four described regulator-set tariffs or recovery on the regulatory asset base without reference to it.  At the time of the survey, Germany \citep{Bundesnetzagentur2024GrandWANDA} and Denmark \citep{KlimaEnergiogForsyningsministeriet2025BrintinfrastrukturSyvtallet.} had proposed capacity tariff schemes based on intertemporal cost allocation, while the other surveyed TSOs indicated that they were undergoing regulatory evaluation of tariff methodologies.

% Acer recommendation
Furthermore, ACER conducted a consultation on hydrogen tariff design \citep{ACER2025ACERInfrastructure}, concluding in July 2025 with an endorsement of intertemporal cost allocation \citep{ACER2025ACERCost-allocation}. Although endorsing the intertemporal cost allocation scheme, ACER stressed that the scheme alone is insufficient, as hydrogen networks lack binding long-term capacity commitments, leaving TSOs without predictable revenues to recover investment costs. Hence, ACER identified three complementary elements. First, national risk-sharing schemes, such as subsidies or pre-booking capacity obligations, are needed to directly reduce volume risk. Second, regulated revenues must reflect only the residual risks borne by TSOs after public support is applied, avoiding double compensation for risks already covered by state instruments. Third, regulatory design must ensure credibility and adaptability through scenario-based demand forecasting, immediate recovery of variable operating costs, and explicit approval of assets subject to intertemporal allocation. The mechanism shifts risk rather than removing it. The regulatory debt accumulated during under-recovery is repaid only if projected demand materialises; if it does not, the balance falls on a smaller user base, producing the tariffs the mechanism was intended to avoid, or on the state. The deferral is long: the German amortisation account is scheduled to balance in 2055 for a network completed in 2032, so tariffs set today would bind users entering the market two decades later. These drawbacks are the counterpart to the three complementary elements identified above, and they explain why endorsing the mechanism alone does not settle the question of who bears the volume risk.

% Congestion management
Regarding congestion management, the expected limited network utilisation in early market phases leads TSOs to anticipate low congestion levels. As a result, congestion management frameworks remain largely under development across surveyed TSOs, and are expected to evolve in parallel with market maturation.
\subsection{Network Codes on Balancing}
The NC balancing survey covered three components: balancing regime design, linepack allocation, and Within-Day Obligations (WDO). Three of the seven respondents had an operative approach to a balancing regime in place at the time of the survey, and, due to the nature of hydrogen and its infrastructure, major changes would be required relative to the natural gas balancing framework. While most TSOs are still developing their approaches, initial designs were available for Denmark, Germany, and the Netherlands.

\begin{figure}
\centering
\includegraphics[width = \linewidth, trim={0cm 10.4cm 0cm 9.8cm},clip]{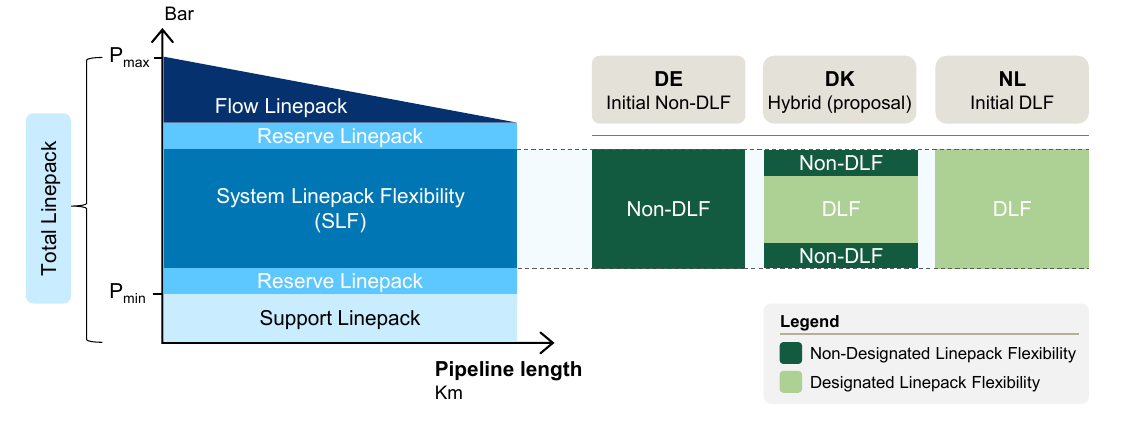}
\caption{Linepack flexibility regions in hydrogen pipeline systems. The diagram illustrates linepack flexibility in the system, subdivided into designated linepack flexibility (DLF), allocated to individual network users, and non-designated linepack flexibility (Non-DLF), available collectively. The pipeline partition illustration is based on \citet{Energinet2024ProposalNetwork}.}
\label{fig:linepack_flexibility_DLF}
\end{figure}

% Balancing regime
The purpose of a balancing regime is to keep pipeline pressure within safe operational bounds. To achieve this, TSOs allocate portions of the total available pipeline pressure to specific purposes. As shown in Figure \ref{fig:linepack_flexibility_DLF}, pipeline pressure is generally organised into four components: support linepack (the minimum cushion gas pressure), reserve linepack held for TSO interventions, system linepack flexibility (SLF) for balancing, and flow linepack providing the pressure differential required for transport. In essence, a balancing regime is defined by how the SLF portion of the pipeline is used for balancing actions. The survey identified three balancing regimes across Germany, Denmark, and the Netherlands, reflected in the Bundesnetzagentur WasABI consultation \citep{Bundesnetzagentur2024SecondWasABi}, the Energinet proposal \citep{Energinet2024ProposalNetwork}, and the Dutch Hynetwork terms and conditions \citep{GasunieNetherlands2026HynetworkConditions}. These balancing regimes differ in their use of a newly introduced approach, namely Designated Linepack Flexibility (DLF), in which the total SLF is partitioned into multiple parts, designated and allocated to individual users. This contrasts with the traditional gas balancing framework, in which SLF is managed at the system level via balancing zones and a helper-causer mechanism. For the purposes of this paper, this traditional approach is referred to as Non-DLF, simply to indicate and differentiate the portion of SLF managed collectively as in gas markets, from the newly introduced DLF approach. The three balancing regimes differ in the allocation of the SLF between DLF and non-DLF, as shown in Figure \ref{fig:linepack_flexibility_DLF}. It is important to note that these designs should be interpreted as transitional. While harmonised hydrogen Network Codes are expected to emerge as the market matures, current approaches reflect early market conditions and are therefore likely to evolve. 

% German proposal
The German WasABI consultation proposes a fully non-DLF balancing regime, aligned with natural gas market design. As in natural gas markets, the system is organised into balancing zones with an accumulated system balance representing the aggregate imbalance of all network users. This value is continuously calculated and published at least every 15 minutes \citep{Bundesnetzagentur2024SecondWasABi}. The framework adopts the familiar green, yellow, and red zone structure. In the green zone, system conditions are stable, and no balancing actions are required. In the yellow zone, the system enters a critical state, and the helper-causer mechanism is activated: network users whose imbalance has the same sign as the system imbalance are classified as causers, while those with an opposite sign are helpers. In the red zone, system conditions become critical and require immediate TSO intervention. This may include the procurement of balancing energy or the use of non-commercial measures. When balancing energy is procured, costs are allocated to the causers based on their contribution to the imbalance. If balancing energy is not procured, penalties are calculated using the EEX HYDRIX index \citep{EEX2026HYDRIX:Hydrogen}, which aims to represent the wholesale hydrogen price in Germany. As the hydrogen market develops, the index is expected to become more representative as market liquidity increases. Penalty revenues are redistributed to helpers.

% Danish proposal
The Danish TSO has proposed a hybrid balancing model combining DLF and Non-DLF elements \citep{Energinet2024ProposalNetwork}. In the hybrid design, total system linepack flexibility consists of a DLF portion and a shared Non-DLF portion. The DLF share is allocated to shippers on a pro rata basis, proportional to booked capacity. For example, a shipper holding 10\% of total capacity receives access to 10\% of the available DLF. Shippers are allowed to maintain an imbalanced portfolio position within this allocated flexibility without incurring penalties. DLF is also expected to be tradable between users, enabling market-based reallocation of flexibility. The remaining SLF forms a shared Non-DLF portion. If multiple users exceed their individual DLF allocations in the same direction, the accumulated system balance approaches the boundary of the total DLF range. Before the reserve linepack is reached and TSO intervention becomes necessary, the Non-DLF portion provides additional buffer at the system level. If the accumulated system balance exceeds the overall SLF limits, balancing actions are triggered, and the users contributing to the imbalance are classified as causers and bear the associated costs. The Danish TSO noted that, at the time of the survey, its final choice between a hybrid approach and a full DLF regime has not yet been made. 

% Netherlands
The Dutch hydrogen network proposes a balancing regime based fully on DLF \citep{GasunieNetherlands2026HynetworkConditions}. As previously mentioned, under this approach, system linepack flexibility is allocated to individual network users, allowing them to maintain imbalances within their designated flexibility range without immediate balancing consequences. The Dutch TSO indicates that this design is intended for the initial phase of network development, when the hydrogen system is not yet fully interconnected and market liquidity remains limited. In later stages, once the network becomes more integrated, the system is expected to transition toward a Non-DLF regime based on balancing zones, similar to those used in natural gas markets. 

In conclusion, the three approaches reflect different trade-offs between individual flexibility and collective buffer capacity. Full DLF maximises individual shipper flexibility but offers no collective buffer before reserve linepack is reached. The hybrid approach adds a shared buffer at the cost of smaller individual allocations and is harder to operate in illiquid markets. Non-DLF is operationally familiar from gas markets but places immediate balancing responsibility on shippers. Hydrogen balancing is not yet regulated at EU level, hence current proposed frameworks, e.g., the full-DLF approach, are member-state designs and may differ in implementation even where similar principles are adopted.

% WDO
Within-day obligations are implemented through frequent system updates and flexible nomination processes. In Germany, 15-minute balance updates provide continuous adjustment signals. In Denmark, within-day balancing is enabled through renomination cycles aligned with 15-minute market intervals. The Dutch proposal does not yet define detailed WDO procedures.

\section{Challenges for the design of future hydrogen markets}
\label{sec:challenges}

\begin{figure}
\centering
\includegraphics[width = \linewidth, trim={2cm 10.5cm 2cm 10cm},clip]{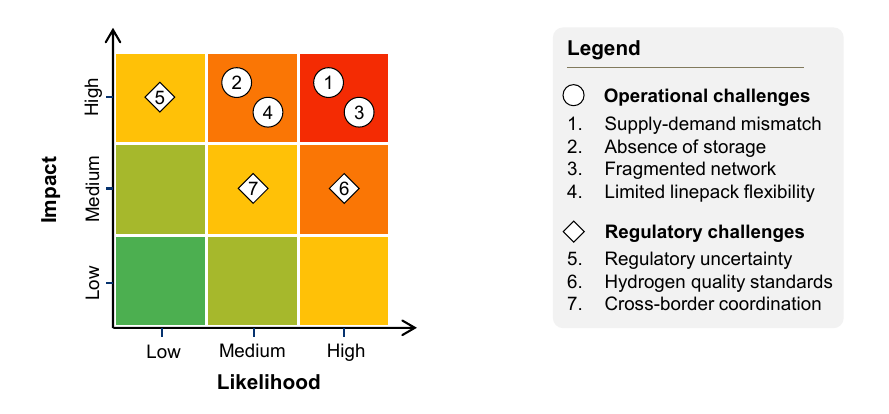}
\caption{Qualitative mapping of hydrogen market challenges by impact and likelihood, assigned by the authors based on the conditions described by participating TSOs.}
\label{fig:risk_matrix}
\end{figure}

The second output of this study identifies seven challenges relevant to early hydrogen market development, grouped into operational and regulatory categories. Operational challenges include a supply-demand mismatch, limited initial storage, a fragmented network, and limited linepack flexibility. Regulatory challenges included regulatory uncertainty, a lack of harmonised hydrogen quality standards, and insufficient cross-border coordination. The challenges were raised by participating TSOs, either in response to direct survey questions or during the broader discussion, and were grouped by the authors into operational and regulatory categories. Impact and likelihood were not elicited as ratings. They were assigned by the authors on a three-level qualitative scale: likelihood, according to whether a condition was described as already present, expected under current plans, or possible but not anticipated; and impact, according to whether the consequence affects individual network users, a national system, or cross-border operations. Figure \ref{fig:risk_matrix} is therefore a qualitative ordering of relative priority, not a quantitative risk assessment.

%Supply-demand mismatch
The most pressing operational challenge is a \textit{temporal mismatch between supply and demand}. In natural gas systems, supply is largely dispatchable, and demand generally follows seasonal cycles. In contrast, with hydrogen this structure is inverted. Consumption from sectors such as steelmaking, chemicals, and fuel production is relatively stable and continuous, while supply is inherently volatile, driven by the intermittent nature of renewable generation rather than operator dispatch. The result is a structural mismatch between supply and demand that could interfere with normal system operations.
% absence of storage
The \textit{absence of large-scale storage} contracted at the TSO-level amplifies this risk. Without it, system-level balancing relies mainly on pipeline linepack. Although ad hoc proposals exist to allocate linepack flexibility among users \citep{Energinet2024ProposalNetwork}, available linepack is expected to remain limited in early market phases. Some TSOs report active development of salt cavern storage, but high capital costs and limited geological availability constrain its deployment.
% Fragmented network
Furthermore, early hydrogen networks are expected to be \textit{fragmented networks} organised into regional clusters, each functioning as an independent balancing zone with limited cross-border interconnection. This results in low liquidity and a small participant base, complicating efficient operation. While the European Hydrogen Backbone aims to interconnect these clusters over time, the transition period poses the challenge of networks operating in isolation.
%Limited linepack flexibility
Finally, \textit{limited linepack flexibility} amplifies the three challenges described above. Linepack capacity depends on both the gas energy content and the pressure range available in the pipeline. Hydrogen has a volumetric energy density approximately one-fourth that of natural gas, meaning the same pipeline volume stores substantially less energy. Combined with the limited pipeline volumes expected in the early stages of hydrogen network development, this leaves little buffer to absorb supply fluctuations, thereby increasing reliance on demand response, dedicated local storage, or in the worst case, curtailment.

%Operating pressure ranges further constrain this buffer: in the Dutch system, the initial range is 30 to 50 bar against a design pressure of 66.2 bar(g) \citep{GasunieNetherlands2026HynetworkConditions}; in the Danish system, limits range from 50 to 90 barg for newly built pipelines and 35 to 58 barg for converted ones \citep{Energinet2025Danish3}. Together, low energy density and limited pressure ranges leave little margin to absorb supply fluctuations, reinforcing the need for demand response, flexible curtailment, or dedicated storage.

% Regulatory challenges
\textit{Regulatory uncertainty} is rated among the most pressing challenges. Market participants must commit to long-term infrastructure investments without finalised Network Codes for hydrogen and without completed EU-level transposition of the Hydrogen and Decarbonised Gas Market Package. This uncertainty directly affects the bankability of hydrogen projects and may defer investment at a critical stage of market development.
% hydrogen standards
A further challenge concerns \textit{hydrogen quality standards}. TSOs require hydrogen to comply with stringent purity specifications. Once hydrogen with a non-compliant composition enters a shared network, its composition cannot be corrected downstream: it propagates through the system and may affect end users and interconnected networks in neighbouring countries. The absence of harmonised cross-border quality standards, therefore, introduces significant operational challenges across operating areas and technical barriers to trade.
% Cross border coordination
Finally, \textit{cross-border coordination} is closely linked to the two preceding challenges. In the absence of agreed Network Codes, TSOs are currently adopting ad hoc measures covering e.g., balancing regimes, nomination procedures and capacity allocation mechanisms. While such arrangements enable early network operation ahead of EU regulation, divergent approaches risk creating friction upon interconnection, when misaligned regulatory timelines and inconsistent operational procedures may impede cross-border flows. EU-level harmonisation of Network Codes, therefore, remains essential to support infrastructure development at the scale required to meet 2030 and 2040 hydrogen targets.

\section{Conclusions and policy implications}
\label{sec:discussion}

\subsection{Conclusions}
% intro
This study examined how the foundational elements of hydrogen market design should be structured, and the challenges they must address, with the explicit aim of informing ENNOH's development of hydrogen Network Codes. Through structured consultation with seven European gas TSOs representing six countries, two outputs were produced: a mapping of how existing gas market design components could be adapted across five regulatory areas (Figure \ref{fig:h2_market_properties}), and an identification of challenges that early hydrogen market design must account for (Figure \ref{fig:risk_matrix}). 

\begin{table}[h]
\centering
\footnotesize
\caption{Staged view of the transition in hydrogen market design, defined by system conditions rather than dates. The table is the authors' synthesis of the maturity-dependent positions reported in Section \ref{sec:transferability}, not a set of positions elicited from respondents; the storage row, and the congestion row beyond Stage 1, extrapolate from the conditions respondents described.}
\label{tab:roadmap}
\begin{tabularx}{\textwidth}{@{}>{\bfseries\RaggedRight\arraybackslash}p{2.1cm}>{\RaggedRight\arraybackslash}X>{\RaggedRight\arraybackslash}X>{\RaggedRight\arraybackslash}X@{}}
\toprule
Component & \textbf{Stage 1} \newline Isolated clusters & \textbf{Stage 2} \newline Interconnection & \textbf{Stage 3} \newline Integrated market \\
\midrule
Shipper roles & Integrated producers and consumers & Commodity shippers enter & Speculative shippers and traders \\
\addlinespace
Trading & Bilateral over-the-counter, long-term contracts & Over-the-counter with emerging exchange trading on virtual points & Exchange reference prices alongside over-the-counter \\
\addlinespace
Capacity products & Yearly and multi-year commitments & Yearly with a reserved share of short-term products & Full product suite as in gas markets \\
\addlinespace
Allocation & First-come first-served or negotiated terms & First-come first-served with a published auction trigger & Auctions \\
\addlinespace
Tariffs & Intertemporal allocation, regulatory debt accumulating & Debt peak and transition to over-recovery & Cost recovery on the regulatory asset base \\
\addlinespace
Congestion & No mechanism; congestion not expected & Congestion management introduced & Gas-equivalent congestion management \\
\addlinespace
Balancing & Flexibility designated to individual users, short balancing periods & Hybrid designated and system-level flexibility & System-level balancing zones \\
\addlinespace
Storage & Linepack only & First underground storage commissioned & Storage-backed seasonal flexibility \\
\bottomrule
\end{tabularx}
\end{table}

Table \ref{tab:roadmap} organises the maturity-dependent positions from Section 3 into three stages defined by system conditions rather than dates: isolated regional clusters with few participants, interconnection with growing liquidity, and an integrated market with utilisation approaching design capacity. The staging is the authors' synthesis, since respondents described conditions rather than sequences; it is an organising device rather than a forecast, and stages are not synchronised across member states. Four conditions drive the transitions: participant numbers for the shipper and trading transitions, physical interconnection for balancing, utilisation for allocation and congestion, and cumulative booked capacity for tariffs. A transitional arrangement introduced without a stated exit condition tends to persist, which is why each stage is defined by the condition that ends it rather than by a date.

% final take
Taken together, these findings indicate that hydrogen market design cannot simply replicate that of natural gas. Selective and targeted adaptation is both necessary and feasible, but success depends on addressing future operational and regulatory challenges before early markets reach operational scale. The timing of these choices is relevant with a historical perspective. European natural gas markets initially evolved under vertically integrated structures, with competitive frameworks emerging only after successive EU liberalisation packages. By contrast, the EU Gas and Hydrogen Package establishes a regulatory framework for hydrogen from day one. As hydrogen markets will therefore be regulated from the start, this study aims to inform these early design decisions.

\subsection{Policy implications}

Three policy messages follow from the consultation. They are limited to points where the evidence is direct, rather than extended into areas such as storage investment or tariff risk-sharing where the survey gives less to stand on.

% balancing
The three balancing regimes described in Section 3.5 differ because the underlying networks differ, in size, in liquidity, and in access to storage. A single balancing regime imposed before interconnection would require some operators to adopt a design unsuited to their network. The more defensible target for EU-level work is therefore the interface between systems rather than the regime within them: how imbalance is measured, what flexibility a shipper is entitled to, and how imbalance is settled once networks are physically connected. Convergence of the regimes themselves can be left until interconnection makes it necessary.

% hydrogen nomination
Hydrogen nomination should be set relative to the electricity market rather than fixed independently of it, since electrolyser output follows the electricity price. Opening nomination once the day-ahead auction has cleared lets a shipper commit against a known price instead of a forecast, and a renomination window after intraday trading lets the hydrogen schedule follow the electricity schedule within the delivery day. No national framework has yet fixed these gate closures, but the direction indicated by respondents, and already reflected in the Danish proposal, points this way.

% Hydrogen quality
Hydrogen quality specifications should be agreed before cross-border flows begin rather than afterwards. Once hydrogen of a non-compliant composition enters a shared network, it cannot be corrected downstream, and divergent national specifications become a technical barrier to trade at the point of interconnection. The same principle applies to the transitional arrangements described elsewhere in this paper: several are appropriate only while the market is small, and each should be adopted together with the condition under which it is expected to end.
\subsection{Limitations and Further Work}
The findings of this study should be interpreted in light of the following limitations.

% We included only TSO perspectives and not the ones of the other stakeholder
First, the consultation covers one stakeholder family, TSOs, whose interests are not neutral with respect to the components assessed: tariff design determines TSO cost recovery, capacity product duration the commitment they can require of users, and the balancing regime the operational risk they bear relative to shippers. Network operators provide the technical input to the Network Codes through ENNOH, but ACER reviews the drafts and the European Commission adopts them, so the positions reported here are input to that process rather than an outcome of it. A more complete assessment would incorporate perspectives from hydrogen producers, industrial consumers, storage operators, traders, electricity market operators, and national regulatory authorities.

%  Low response rate
Second, the participating operators are a subset of the ENTSOG membership, and Southern and South-Eastern Europe are unrepresented, with the geographic concentration following from non-response rather than from the selection criterion. The sample nonetheless spans a range of market maturity, and the divergences reported in Section 3 track that gradient. If more advanced operators propose greater departures from gas arrangements than less advanced ones, a sample weighted towards them overstates the adaptation required, and the finding that most components transfer with minor or no modification should be read as conservative.

% nice analysis but we don't have a pipeline yet
Finally, the analysis necessarily focuses on market design mechanisms for a hydrogen network that remains largely undeveloped. The findings therefore rely on assumptions about how infrastructure, liquidity, and regulatory frameworks will evolve. As hydrogen pipelines are commissioned and early trading activity emerges, these design choices will become empirically testable and can be refined accordingly.

% Acknowledgments
\newpage

% Author Contributions

\section*{CRediT authorship contribution statement}
\textbf{Marco Saretta}: Conceptualization, Methodology, Investigation, Data curation, Formal analysis, Visualization, Software, Project administration, Writing - original draft. 
\textbf{Enrica Raheli}: Supervision, Writing - review and editing.
\textbf{Jalal Kazempour}: Supervision, Writing - review and editing. 

% Declaration of Interests
\section*{Declaration of competing interest}
The authors declare that they have no known competing financial interests or personal relationships that could have influenced the work reported in this paper.

\section*{Acknowledgments}
This work is conducted as part of the Industrial PhD programme ``Commoditizing green hydrogen in Europe: from efficient market design to optimal contracting and investment'', funded by Innovation Fund Denmark with grant number 4297-00039B, Ramboll Fonden, and Ramboll Denmark A/S. The authors express their gratitude to all EU gas TSOs that participated in this initiative and the survey. Special thanks to Christian Rutherford from the Danish TSO for his thorough discussions, constructive feedback, and final review of the manuscript. Furthermore, the authors acknowledge the feedback provided by several organisers and participants at the Loyola Autumn Research School 2025, organised by the Florence School of Regulation, which contributed to reviewing and enhancing this paper, significantly refining the policy messages. Sincere thanks also go to Alexandra Lüth and Andrea Saretta from Ramboll Denmark for their feedback and help with reviewing the draft; to Miguel Martinez Rodriguez from ACER for his significant contributions to the review of this paper, and for providing valuable feedback on this work; and to Mirco Dian for his help with software curation and data collection.

\section*{Data availability}
The survey data that support the findings of this study are available from the corresponding author upon request.
Code and supplementary material are publicly available at \cite{Saretta2026DesigningRepository}.

\section*{Declaration of generative AI and AI-assisted technologies in the writing process}
During the preparation of this work the authors used generative AI models from OpenAI and Anthropic for translation of original sources and for grammatical revisions. After using these tool, the authors thoroughly reviewed and edited all outputs and take full responsibility for the content of this paper.

% References
\newpage
\bibliographystyle{elsarticle-harv} 
\bibliography{references}

@techreport{ACER2025ACERInfrastructure,
    title = {{ACER Public consultation on inter-temporal cost allocation mechanisms for financing hydrogen infrastructure}},
    year = {2025},
    author = {{ACER}},
    month = {7},
    url = {https://www.acer.europa.eu/documents/public-consultations/pc2025g01-public-consultation-inter-temporal-cost-allocation-mechanisms-financing-hydrogen-infrastructure}
}

@techreport{ACER2025ACERCost-allocation,
    title = {{ACER Recommendation on inter-temporal cost-allocation}},
    year = {2025},
    author = {{ACER}},
    month = {7}
}

@techreport{ENTSOG2018BalancingOverview,
    title = {{Balancing Network Code – An Overview}},
    year = {2018},
    author = {{ENTSOG}},
    month = {3}
}

@article{Lagioia2023BlueAffairs,
    title = {{Blue and green hydrogen energy to meet European Union decarbonisation objectives. An overview of perspectives and the current state of affairs}},
    year = {2023},
    journal = {International Journal of Hydrogen Energy},
    author = {Lagioia, Giovanni and Spinelli, Maria Pia and Amicarelli, Vera},
    number = {4},
    month = {1},
    pages = {1304--1322},
    volume = {48},
    publisher = {Pergamon},
    url = {https://www.sciencedirect.com/science/article/abs/pii/S0360319922046675},
    doi = {10.1016/J.IJHYDENE.2022.10.044},
    issn = {0360-3199}
}

@techreport{KlimaEnergiogForsyningsministeriet2025BrintinfrastrukturSyvtallet.,
    title = {{Brintinfrastruktur til tyskland: muligg{\o}relse af syvtallet.}},
    year = {2025},
    author = {{Klima Energi og Forsyningsministeriet}},
    month = {2},
    url = {https://www.kefm.dk/Media/638744580036601360/Aftale%20om%20brintinfrastruktur%20til%20Tyskland_muligg%C3%B8relse%20af%20Syvtallet.pdf}
}

@techreport{ENTSOG2017CapacityCode,
    title = {{Capacity Allocation Mechanisms Network Code}},
    year = {2017},
    author = {{ENTSOG}}
}

@misc{Energinet2025CapacityBackbone,
    title = {{Capacity sale for the Danish Hydrogen Backbone}},
    year = {2025},
    author = {{Energinet}},
    url = {https://en.energinet.dk/hydrogen/for-hydrogen-market-players/capacity-sale-hydrogen/}
}

@techreport{Fraile2025Clean2025,
    title = {{Clean Hydrogen Monitor 2025}},
    year = {2025},
    author = {Fraile, Daniel and Muron, Matus and Pawelec, Grzegorz and Santos, Sara and Staudenmayer, Olivia},
    url = {https://hydrogeneurope.eu/wp-content/uploads/2025/09/Clean_Hydrogen_Monitor_09-2025_DIGITAL.pdf},
    institution = {Hydrogen Europe}
}

@article{Guinlle2026CompetingMarket,
    title = {{Competing corporate visions, ambiguous policy proposals: Complexities surrounding the emerging EU hydrogen market}},
    year = {2026},
    journal = {Energy Policy},
    author = {Guinlle, Gabriella and Kolk, Ans},
    month = {11},
    publisher = {Elsevier Ltd},
    doi = {10.1016/j.enpol.2026.115463},
    issn = {03014215}
}

@article{Nunez-Jimenez2022CompetitiveUnion,
    title = {{Competitive and secure renewable hydrogen markets: Three strategic scenarios for the European Union}},
    year = {2022},
    journal = {International Journal of Hydrogen Energy},
    author = {Nu{\~{n}}ez-Jimenez, Alejandro and De Blasio, Nicola},
    number = {84},
    month = {10},
    pages = {35553--35570},
    volume = {47},
    publisher = {Elsevier Ltd},
    doi = {10.1016/j.ijhydene.2022.08.170},
    issn = {03603199}
}

@article{Brandt2024CostFramework,
    title = {{Cost and competitiveness of green hydrogen and the effects of the European Union regulatory framework}},
    year = {2024},
    journal = {Nature Energy},
    author = {Brandt, Jonathan and Iversen, Thore and Eckert, Christoph and Peterssen, Florian and Bensmann, Boris and Bensmann, Astrid and Beer, Michael and Weyer, Hartmut and Hanke-Rauschenbach, Richard},
    number = {6},
    month = {6},
    pages = {703--713},
    volume = {9},
    publisher = {Nature Research},
    doi = {10.1038/s41560-024-01511-z},
    issn = {20587546}
}

@techreport{EuropeanParliament2023DelegatedOrigin,
    title = {{Delegated Regulation (EU) 2023/1184 of 10 February 2023 supplementing Directive (EU) 2018/2001 of the European Parliament and of the Council by establishing a Union methodology setting out detailed rules for the production of renewable liquid and gaseous transport fuels of non-biological origin}},
    year = {2023},
    author = {{European Parliament}},
    month = {2},
    url = {https://eur-lex.europa.eu/eli/reg_del/2023/1184/oj/eng}
}

@article{Hasan2025DesigningReform,
    title = {{Designing effective hydrogen Markets: Policy recommendations from electricity and gas market reform}},
    year = {2025},
    journal = {Energy Policy},
    author = {Hasan, Shahid and Fuentes, Rolando},
    month = {12},
    volume = {207},
    publisher = {Elsevier Ltd},
    doi = {10.1016/j.enpol.2025.114816},
    issn = {03014215}
}

@techreport{Saretta2026DesigningRepository,
    title = {{Designing European Hydrogen Markets. GitHub repository}},
    year = {2026},
    author = {{Saretta}},
    url = {https://github.com/marco-saretta/hydrogen-market-design}
}

@misc{YahooFinance2025DutchData,
    title = {{Dutch TTF Natural Gas Calendar, Historical Prices {\&} Data}},
    year = {2025},
    author = {{Yahoo Finance}},
    url = {https://finance.yahoo.com/quote/TTF%3DF/history/}
}

@misc{ENTSOGENTSOGPlatform,
    title = {{ENTSOG - Transparency Platform}},
    author = {{ENTSOG}},
    url = {https://transparency.entsog.eu}
}

@misc{EuropeanCommission2024EUPackage,
    title = {{EU Hydrogen and Decarbonized Gas Market Package}},
    year = {2024},
    author = {{European Commission}},
    month = {5},
    url = {https://energy.ec.europa.eu/topics/markets-and-consumers/hydrogen-and-decarbonised-gas-market_en}
}

@misc{EuropeanHydrogenObservatory2025EUDeal,
    title = {{EU Hydrogen Strategy under the EU Green Deal}},
    year = {2025},
    author = {{European Hydrogen Observatory}},
    url = {https://observatory.clean-hydrogen.europa.eu/eu-policy/eu-hydrogen-strategy-under-eu-green-deal#:~:text=The%20target%20is%2040%20GW}
}

@techreport{EuropeanHydrogenBackbone2024EuropeanCompetitiveness,
    title = {{European Hydrogen Backbone: Boosting EU Resilience and Competitiveness. European Hydrogen Backbone: Boosting EU Resilience and Competitiveness}},
    year = {2024},
    author = {{European Hydrogen Backbone}}
}

@techreport{ACER2025EuropeanReport,
    title = {{European hydrogen markets - 2025 Monitoring Report}},
    year = {2025},
    author = {{ACER}},
    month = {12},
    url = {https://www.acer.europa.eu/},
    institution = {ACER}
}

@article{Bundesnetzagentur2024GrandWANDA,
    title = {{Grand Ruling Chamber for Energy on decision [GBK-24-01-2{\#}1] (WANDA)}},
    year = {2024},
    author = {{Bundesnetzagentur}},
    month = {4}
}

@misc{EEX2026HYDRIX:Hydrogen,
    title = {{HYDRIX: First market-based index for hydrogen}},
    year = {2026},
    author = {{EEX}},
    url = {https://www.eex.com/en/markets/hydrogen}
}

@misc{EuropeanCommission2025HydrogenMarket,
    title = {{Hydrogen and decarbonised gas market}},
    year = {2025},
    author = {{European Commission}},
    url = {https://energy.ec.europa.eu/topics/markets-and-consumers/hydrogen-and-decarbonised-gas-market_en}
}

@misc{ENTSOG2026HydrogenMap,
    title = {{Hydrogen Infrastructure Map}},
    year = {2026},
    author = {{ENTSOG} and {GIE} and {EUROGAS} and {GEODE} and {GD4S} and {CEDEC}},
    url = {https://www.h2inframap.eu/#map}
}

@misc{EuropeanHydrogenObservatory2026HydrogenProduction,
    title = {{Hydrogen production}},
    year = {2026},
    author = {{European Hydrogen Observatory}},
    url = {https://observatory.clean-hydrogen.europa.eu/hydrogen-landscape/production-trade-and-cost/hydrogen-production}
}

@misc{GasunieNetherlands2026HynetworkConditions,
    title = {{Hynetwork Terms and Conditions}},
    year = {2026},
    author = {{Gasunie Netherlands}},
    url = {https://www.hynetwork.nl/zakelijk/klant-worden/contracten}
}

@article{Sargent2025LinepackModelling,
    title = {{Linepack flexibility for hydrogen and natural gas pipes: A new flexibility metric for energy systems modelling}},
    year = {2025},
    journal = {International Journal of Hydrogen Energy},
    author = {Sargent, Philip and Sargent, Michael},
    month = {5},
    pages = {139--142},
    volume = {132},
    publisher = {Pergamon},
    url = {https://www.sciencedirect.com/science/article/abs/pii/S0360319925018543},
    doi = {10.1016/J.IJHYDENE.2025.04.197},
    issn = {0360-3199}
}

@article{Niedrig2024MarketMarket,
    title = {{Market Design Options for a Hydrogen Market}},
    year = {2024},
    author = {Niedrig, Nicolas and Giehl, Johannes and Jahnke, Philipp and M{\"{u}}ller-Kirchenbauer, Joachim}
}

@article{Martinez-Rodriguez2026PipelineNetworks,
    title = {{Pipeline regulation for hydrogen: choosing between paths and networks}},
    year = {2026},
    journal = {Energy Policy},
    author = {Martinez-Rodriguez, Miguel and Chyong, Chi Kong and Fitzgerald, Timothy and Vazquez, Miguel and Hidalgo, Antonio},
    month = {1},
    pages = {114846},
    volume = {208},
    publisher = {Elsevier},
    url = {https://www.sciencedirect.com/science/article/pii/S0301421525003532?ref=pdf_download&fr=RR-2&rr=9ba3ed790885ed4f},
    doi = {10.1016/J.ENPOL.2025.114846},
    issn = {0301-4215}
}

@misc{PRISMA2025PRISMAPlatform,
    title = {{PRISMA gas capacity trading platform}},
    year = {2025},
    author = {{PRISMA}},
    url = {https://www.prisma-capacity.eu/}
}

@techreport{ENTSOE2023ProcessGuidelines,
    title = {{Process for developing network codes and amendments to network codes and  guidelines}},
    year = {2023},
    author = {{ENTSOE}},
    url = {www.entsoe.eu}
}

@techreport{Energinet2024ProposalNetwork,
    title = {{Proposal on balancing the Danish hydrogen transmission network}},
    year = {2024},
    author = {{Energinet}},
    month = {6},
    url = {https://energinet.dk/media/a4gdh34h/proposal-on-balancing-the-danish-hydrogen-transmission-network.pdf}
}

@article{Johnson2025RealisticTransition,
    title = {{Realistic roles for hydrogen in the future energy transition}},
    year = {2025},
    journal = {Nature Reviews Clean Technology 2025 1:5},
    author = {Johnson, Nathan and Liebreich, Michael and Kammen, Daniel M. and Ekins, Paul and McKenna, Russell and Staffell, Iain},
    number = {5},
    month = {4},
    pages = {351--371},
    volume = {1},
    publisher = {Nature Publishing Group},
    url = {https://www.nature.com/articles/s44359-025-00050-4},
    doi = {10.1038/s44359-025-00050-4},
    issn = {3005-0685}
}

@misc{EuropeanParliament2024RegulationHydrogen,
    title = {{Regulation on the internal markets for renewable gas, natural gas and hydrogen, }},
    year = {2024},
    author = {{European Parliament}},
    url = {https://eur-lex.europa.eu/legal-content/EN/TXT/?uri=CELEX%3A32024R1789}
}

@misc{GasunieTransportServices2025SBSPOS,
    title = {{SBS and POS}},
    year = {2025},
    author = {{Gasunie Transport Services}},
    url = {https://www.gasunietransportservices.nl/en/shipper-trader/balancing-regime/sbs-and-pos}
}

@article{Bundesnetzagentur2024SecondWasABi,
    title = {{Second consultation in the determination proceedings for a basic compensation and balancing model for hydrogen (“WasABi”)}},
    year = {2024},
    author = {{Bundesnetzagentur}},
    month = {12}
}

@article{Bundesnetzagentur2024SecondWaKandA,
    title = {{Second consultation in the determination proceedings on a basic model for hydrogen capacity and managing network access ("WaKandA")}},
    year = {2024},
    author = {{Bundesnetzagentur}},
    month = {12}
}

@techreport{ENTSOG2018TariffOverview,
    title = {{Tariff Network Code – An Overview}},
    year = {2018},
    author = {{ENTSOG}},
    month = {9}
}

@techreport{EuropeanHydrogenObservatory2025TheLandscape,
    title = {{The European hydrogen market landscape}},
    year = {2025},
    author = {{European Hydrogen Observatory}},
    month = {11}
}

@article{Steinbach2024TheTrading,
    title = {{The future European hydrogen market: Market design and policy recommendations to support market development and commodity trading}},
    year = {2024},
    journal = {International Journal of Hydrogen Energy},
    author = {Steinbach, Sarah A. and Bunk, Nikolas},
    month = {6},
    pages = {29--38},
    volume = {70},
    publisher = {Elsevier Ltd},
    doi = {10.1016/j.ijhydene.2024.05.107},
    issn = {03603199}
}

@article{Neumann2023TheEurope,
    title = {{The potential role of a hydrogen network in Europe}},
    year = {2023},
    journal = {Joule},
    author = {Neumann, Fabian and Zeyen, Elisabeth and Victoria, Marta and Brown, Tom},
    number = {8},
    month = {8},
    pages = {1793--1817},
    volume = {7},
    publisher = {Cell Press},
    url = {https://www.cell.com/action/showFullText?pii=S2542435123002660 https://www.cell.com/action/showAbstract?pii=S2542435123002660 https://www.cell.com/joule/abstract/S2542-4351(23)00266-0},
    doi = {10.1016/J.JOULE.2023.06.016/ASSET/66A7C216-C041-4A9B-B0CC-DB2B9414FD73/MAIN.ASSETS/GR5.JPG},
    issn = {25424351},
    arxivId = {2207.05816}
}

@article{Shen2024TheImpacts,
    title = {{The role of hydrogen in iron and steel production: Development trends, decarbonization potentials, and economic impacts}},
    year = {2024},
    journal = {International Journal of Hydrogen Energy},
    author = {Shen, Jialin and Zhang, Qi and Tian, Shuoshuo and Li, Xingyu and Liu, Juan and Tian, Jinglei},
    month = {11},
    pages = {1409--1422},
    volume = {92},
    publisher = {Pergamon},
    url = {https://www.sciencedirect.com/science/article/abs/pii/S0360319924045798},
    doi = {10.1016/J.IJHYDENE.2024.10.368},
    issn = {0360-3199}
}

@misc{EuropeanCommission2009ThirdPackage,
    title = {{Third Energy Package}},
    year = {2009},
    author = {{European Commission}},
    url = {https://energy.ec.europa.eu/topics/markets-and-consumers/governance-internal-energy-market_en}
}

% Appendix (with continuous arabic page numbering)
\newpage
\appendix
\newpage
\section{TSO Hydrogen Survey}
\label{sec:survey}

\subsection{Survey administration}

The 37-question survey was structured around five components derived from EU gas market regulation: (A) market and regulatory design principles from Regulation 2024/1789, (B) trading framework, (C) Network Code on Capacity Allocation Mechanisms (NC CAM), (D) Network Code on Tariffs (NC TAR), and (E) Network Code on Balancing (NC BAL). Questions used a mix of open-ended and yes/no formats to capture detailed qualitative reasoning. The survey was pilot-tested with the Danish TSO under agreement, in August 2024 and administered between September 2024 and February 2025.

The ENTSOG directory distinguishes three categories of participant: members, associated partners, and observers. The survey was addressed to the 43 members, on the criterion set out in Section 2.2. Each member was approached through a named contact person where one was available, otherwise through the contact form or information address published on the operator's website, or otherwise through the operator's general corporate address. Members were contacted between September 2024 and February 2025.

The survey was distributed in Excel format with an estimated completion time of 45–90 minutes. Respondents returned completed questionnaires via email. Where responses required clarification or additional context, follow-up discussions were conducted through video call on Microsoft Teams to complement the written answers.

Seven TSOs from six countries provided responses: Belgium, Denmark, France, Germany, Lithuania, and the Netherlands. Participating countries account for around 45\% of the planned European hydrogen pipeline length by 2040 according to European Hydrogen Backbone data \citep{ENTSOG2026HydrogenMap}. Respondents held positions in operations management, strategic planning, and regulatory affairs. Anonymity was granted both to the TSO organisations and individual participants, who were anonymised, with results reported only at the country level.

\subsection{Survey questions}

\begin{table}[h!]
\centering
\caption{Complete survey distributed for TSO consultations on hydrogen market design.}
\label{tab:survey_questions}
\footnotesize
\setlength{\tabcolsep}{4pt}
\begin{tabularx}{\linewidth}{>{\raggedright\arraybackslash}p{0.8cm} >{\raggedright\arraybackslash}X}
\toprule
\textbf{ID} & \textbf{Question} \\
\midrule
\multicolumn{2}{l}{\textbf{A. Market and Regulatory Design (13 questions)}} \\
\midrule
A.1 & Who is the regulating authority in your country? \\
A.2 & Do you have a regulation scheme in place? If so, when is it expected to be implemented? \\
A.3 & Will the entry/exit point scheme be maintained in the hydrogen market? \\
A.4 & Will the gas day and gas year structure be kept for hydrogen? \\
A.5 & How will the shipper roles evolve, especially at early market stages with lower liquidity? \\
A.6 & In your country, how would you expect the hydrogen generation to be? Green, blue, or gray? \\
A.7 & How will Guarantees of Origin play a role with RFNBO requirements? \\
A.8 & What technical and regulatory challenges of the new hydrogen economy are specific to your TSO's country? \\
A.9 & What is the hydrogen strategy for the new hydrogen economy in your TSO's country? \\
A.10 & What are the expectations regarding interaction with TSOs within the same country and neighbouring countries? \\
A.11 & Which role will DSOs play at the early market stage, given the large industrial nature of \ce{H2} off-takers? \\
A.12 & What in your opinion has the most urgent requirement to be developed? \\
A.13 & If you had time to do research, what would be some challenges that you wish you could solve for future hydrogen markets?\\

\midrule
\multicolumn{2}{l}{\textbf{B. Trading Framework (7 questions)}} \\
\midrule
B.1 & Should we expect an exchange for \ce{H2} (like ICE)? \\
B.2 & Should we expect a virtual trading hub for \ce{H2} (like TTF)? \\
B.3 & Should we expect \ce{H2} to be traded at exchanges, via OTC trades, or both, as in natural gas? \\
B.4 & How will the price formation mechanism function? \\
B.5 & How correlated will hydrogen prices and electricity prices be? \\
B.6 & Which role will derivative contracts (futures, forwards, PPAs) play in electricity provision and hydrogen intake? \\
B.7 & What will be the impact of \ce{CO2} quotas? \\
\midrule
\multicolumn{2}{l}{\textbf{C. Network Code on Capacity Allocation Mechanisms (NC CAM) (7 questions)}} \\
\midrule
C.1 & Will the scheme of nominations, re-nomination, and allocations be kept for hydrogen networks? \\
C.2 & How will the gate closure time for nomination submissions in the gas day be affected? \\
C.3 & Will the capacity products differ from natural gas markets (e.g., yearly, quarterly, monthly, daily, intraday)? \\
C.4 & What is the reasoning behind the capacity products you offer? \\
C.5 & How much capacity is expected to be available in the early stages? \\
C.6 & How and when are you scheduled to go live with the capacity? \\
C.7 & Will the TSO assign capacity with mechanisms different from ascending clock auctions? \\
\midrule
\multicolumn{2}{l}{\textbf{D. Network Code on Tariffs (NC TAR) (2 questions)}} \\
\midrule
D.1 & How will the capacity base price and premium be decided? (System tariffs) \\
D.2 & Will the current congestion management system need adjustments for hydrogen? \\
\midrule
\multicolumn{2}{l}{\textbf{E. Network Code on Balancing (NC BAL) (8 questions)}} \\
\midrule
E.1 & Do you have a network code on balancing which will govern the framework for the balancing model? \\
E.2 & How will the hydrogen target balancing model differ from the natural gas one? \\
E.3 & What is the expected role of linepack in hydrogen markets? Can you quantify it? \\
E.4 & Will the determination methodology for green, yellow, and red zones be adjusted in the future \ce{H2} market? \\
E.5 & How will continuous balancing work for hydrogen markets? \\
E.6 & How can the duality between storage needs and the lack of availability at early market stages be addressed? \\
E.7 & How will shippers be incentivized to maintain an optimal balance in the network? \\
E.8 & How will shippers deal with variable supply and constant offtake? \\
\bottomrule
\end{tabularx}
\end{table}

\newpage
\section{EU Gas Market Principles in the EU}
\label{sec:gas_market}

This section provides an extensive description of the foundational mechanisms of the EU gas market. Given the absence of comprehensive literature on the holistic functioning of gas markets, this section serves two purposes: centralising market operations and identifying the key market design components for gaseous commodities in the EU. The analysis covers market structure, capacity allocation, nomination procedures, balancing frameworks, and commodity trading to establish five core market design components that will guide the hydrogen market analysis.

\subsection{Entry-exit structure}
The natural gas system is intended for the exchange of natural gas, achieved via a network of pipelines operated by the TSO. TSOs are responsible for the safe and reliable operation and maintenance of the gas network within a designated area defined as a \textit{balancing zone}. Notably, a single balancing zone can be managed by multiple TSOs; e.g., Germany's balancing zone, which aligns with the national borders, includes several TSOs operating across the country. Historically, balancing zones have aligned with national borders. However, the recent introduction of joint balancing zones, spanning multiple national borders, has changed this. Examples include the Danish-Swedish joint balancing zone established in 2019, as well as the joint balancing zones between Belgium and Luxembourg and among Estonia, Latvia, and Finland. While the TSO is responsible for the high-pressure infrastructure, its operation is coupled to the Distribution System Operator (DSO), which operates the low-pressure distribution for end use and delivery to consumers. 

The natural gas system is structured as an \textit{entry-exit system}, in which gas enters or exits the network at specific points designated by the TSO for the respective balancing zone. These points are classified as either entry points, where gas is injected into the system, or exit points, where gas is withdrawn, with some points serving both. Entry and exit points can be located within the balancing zone intended for internal gas exchange, or at its borders, facilitating cross-border exchange among neighbouring zones. Figure \ref{fig:entry-exit} illustrates the joint balancing zone of Denmark and Sweden, displaying only the cross-border entry/exit points, while omitting internal ones.

\begin{figure}
    \centering
    \includegraphics[width=0.9\linewidth]{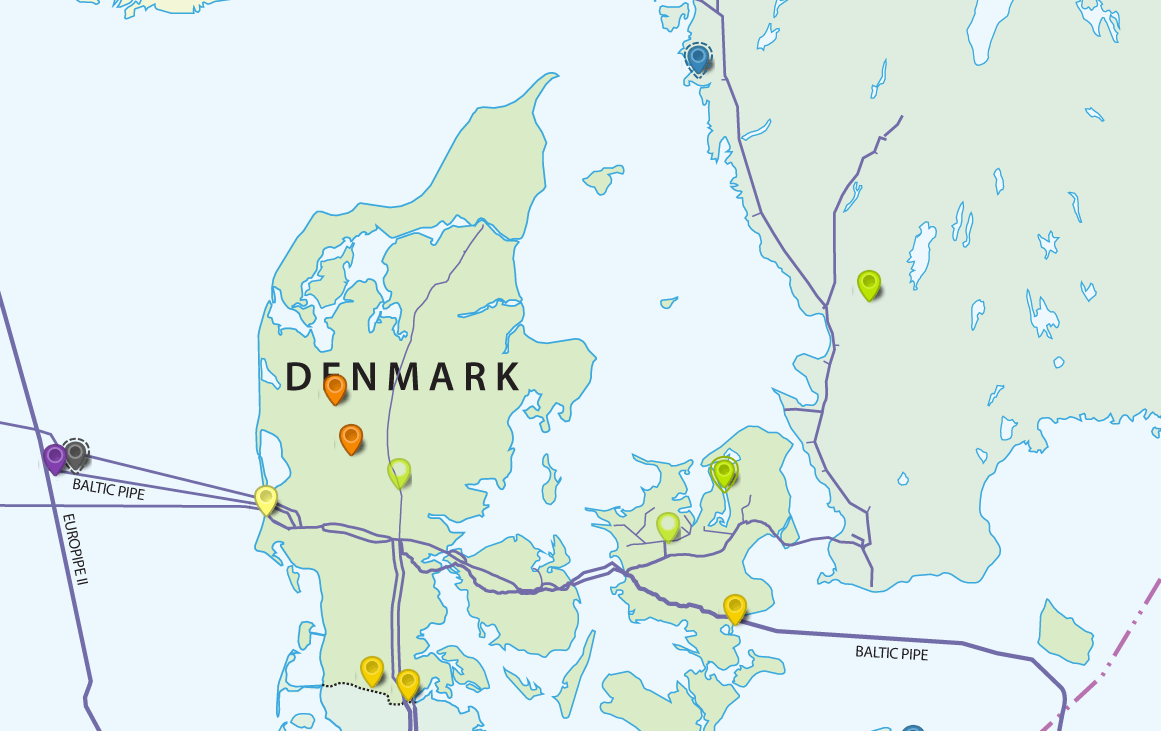}
    \caption{\textbf{The joint balancing zone of Denmark and Sweden.} Gas pipelines are shown in purple, storage points that serve as entry/exit points are in green, the virtual trading point is marked in orange, and cross-border entry/exit points are highlighted in yellow. Entry/exit points within the balancing zone are not displayed.}
    \label{fig:entry-exit}
\end{figure}

In the natural gas system, \textit{network users}, also known as shippers, are parties who transport and exchange natural gas through the network. These include producers, who inject gas into the system; consumers, who withdraw gas at consumption points; and traders, who exchange gas within the market to leverage arbitrage opportunities. 

\subsection{Capacity} 

Each entry and exit point $u$ within the natural gas system set $U$ is characterised by capacity $q_u^{\text{TC}}$, which is the maximum flow of natural gas that can be injected or withdrawn at that point, measured in MWh/h. Users must hold a capacity contract to obtain access rights for any network point, specifying the booked capacity and its duration. For example, a yearly capacity contract of 10 GWh/h at an exit point allows a user to withdraw up to 10 GWh in any given hour during that year. Holding capacity grants the right, but not the obligation, to flow gas at the specified point. Capacity products are available for purchase through online platforms, with PRISMA being the principal platform at the EU level \citep{PRISMA2025PRISMAPlatform}. Notably, PRISMA is not the sole platform available; Poland and Hungary have also established official EU capacity platforms. These capacity products vary in duration (yearly, quarterly, monthly, daily, within-day) and type (firm or interruptible) and are regulated by the EU network codes \citep{Entsog2017CapacityCode}.

\textit{Firm capacity} products guarantee uninterrupted access to the contracted volume for the entire duration of the agreement. Firm capacity contracts are auctioned on PRISMA. The total technical capacity at a network point $u$ is auctioned, and the capacity purchased and booked by shipper $i$ is defined as firm booked capacity $q_{u, i}^{\text{FB}}$. In contrast, the remaining unbooked portion is defined as the firm's available capacity, $q_u^{\text{FA}}$. For any point $u$, shipper $i \in \mathcal{P}$, these capacities are related by:

\begin{equation}
    q_u^{\text{TC}} = q_u^{\text{FA}} + \sum_{i \in \mathcal{P}_u} q_{i, u}^{\text{FB}}, \quad \forall u \in U,\,  i \in \mathcal{P}_u
    \label{eq:firm_capacity}
\end{equation}

where $\mathcal{P}_u \subseteq \mathcal{P}$ is the set of all shippers holding firm capacity contracts at point $u$.

\begin{figure}
    \centering
    \includegraphics[width=\linewidth]{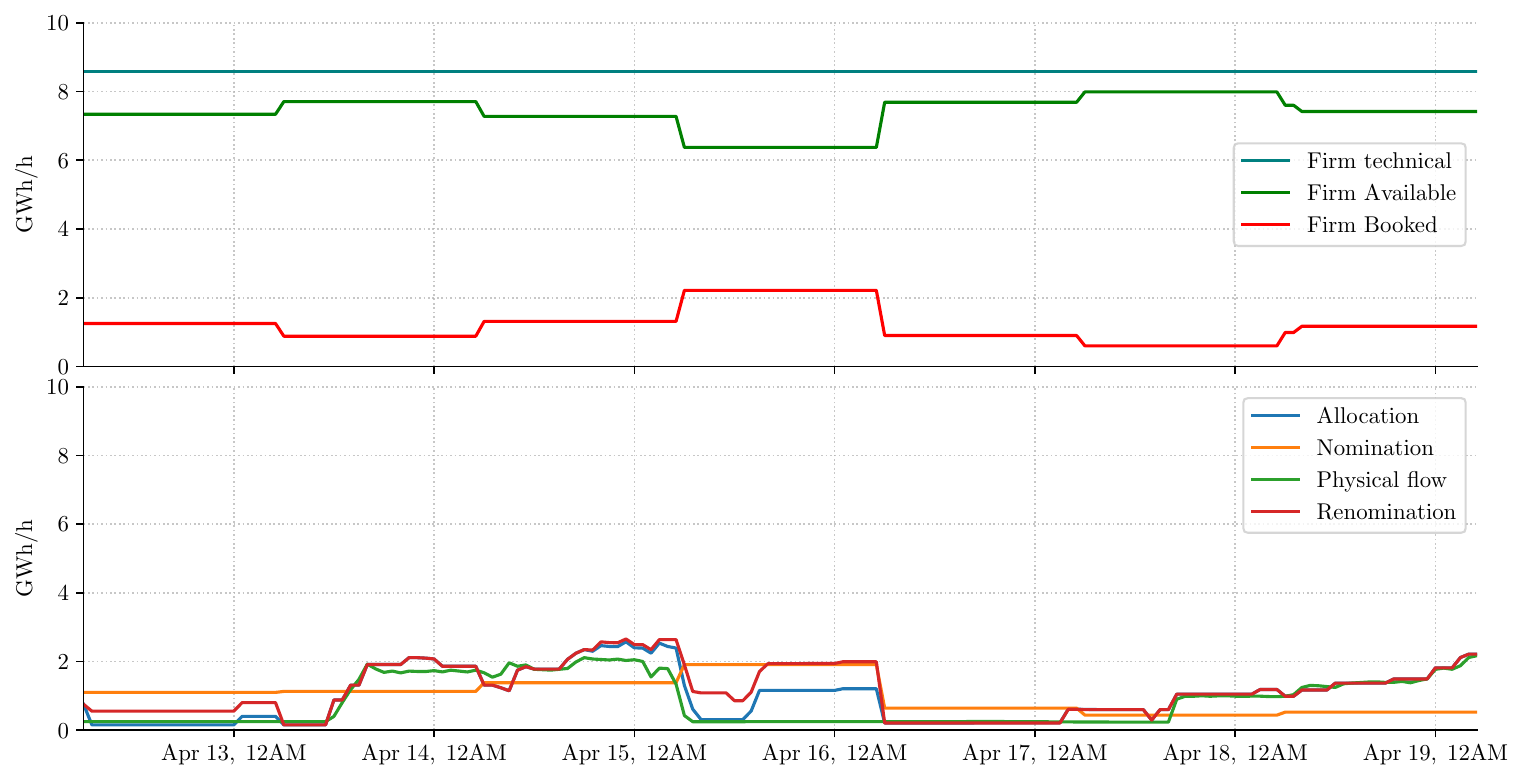}
    \caption{\textbf{Main figures for the Danish entry point of Nybro}. Data displayed from 12 April 2025 06:00 to 19 April 2025 05:00 in hourly resolution \cite{ENTSOGENTSOGPlatform}. The upper subplot shows the breakdown of firm capacity, whereas the lower subplot displays scheduling and allocation.}
    \label{fig:capacity-allocation}
\end{figure}
In Figure \ref{fig:capacity-allocation}, the upper plot illustrates historical capacity data for the Nybro entry point. The technical capacity remains constant at 8,589 GWh, while booked and available capacity fluctuate over time. The sum remains equal to the technical capacity.

\textit{Interruptible capacity} is offered by TSOs to monetise or reschedule availability. These products grant conditional access subject to curtailment when firm capacity is fully utilised, allowing the firm to utilise other infrastructure while preserving the right to interrupt service during congestion. Interruptible capacity is determined through three methods: first, as unused firm capacity; second, via counterflow-based assessment, where TSOs analyse historical flows in the opposite direction to estimate residual capacity; and lastly, based on statistical scenarios, which calculate the likelihood of interruption according to expected flow conditions within the system. The total interruptible capacity $q_u^{\text{IT}}$ at node $u$ is defined as the sum of interruptible capacity booked by shippers $q_{i,u}^{\text{IB}}$ and the remaining available interruptible capacity $q_u^{\text{IA}}$. These capacities are related as follows:

\begin{equation}
    q_u^{\text{IT}} = q_u^{\text{IA}} + \sum_{i \in \mathcal{P}_u} q_{i,u}^{\text{IB}} , \quad \forall u \in U, \, i \in \mathcal{P}_u
\end{equation}

For instance, if the technical capacity $q_u^{\text{TC}}$ is 10 GWh/h, firm booked capacity $q_u^{\text{FB}}$ is 10 GWh/h. Interruptible booked capacity $q_u^{\text{IB}}$ is 5 GWh/h, firm users flowing 5 GWh and interruptible users flowing 5 GWh stay within system limits. However, if firm usage increases to 10 GWh, the TSO curtails interruptible flows to respect the technical capacity limit, leading shippers to interrupt their flows.

\textit{Capacity payment} is defined as the product of the capacity booked by the shipper $i$ at node $u$ and its booking price $\pi_u^{\text{book}}$, applicable at each entry or exit node $u$. When a user $i$ purchases a capacity amount $q_{i,u}$ at node $u$, the resulting capacity payment is calculated as follows:

\begin{equation}
    C_{i,u} = q_{i,u} \pi_{u}^{\text{book}}
    \label{eq:capacity_fee}
\end{equation}

Capacity price, also referred to as capacity tariff, is regulated under the EU Tariff Codes \citep{Entsog2018TariffOverview}, which establish general principles and set upper bounds for short-term pricing multipliers; individual TSOs determine the specific tariffs. To incentivise loyalty from shippers, TSOs commonly differentiate prices based on contract duration. For example, in Denmark and Germany, annual contracts are priced at the base rate $\pi_{u}^{\text{book}}$, whereas shorter-term contracts include premiums: quarterly at $1.1\pi_{u}^{\text{book}}$, monthly at $1.25\pi_{u}^{\text{book}}$, and daily at $1.4\pi_{u}^{\text{book}}$. These figures are not universally applied across the EU but are subject to national implementation within the limits defined by the EU framework. These capacity payments generally account for less than 10\% of the total cost of gas delivery, with approximately 90\% attributable to the cost of the gas itself, including taxes. Capacity payments constitute the primary revenue stream for TSOs, facilitating the operation and maintenance of the gas network. Revenue regulation is managed through the Regulatory Asset Base (RAB) framework, designed to ensure full cost recovery. If actual revenues exceed the allowable level (over-recovery), future tariffs are correspondingly reduced. Conversely, if revenues fall below the allowable level (under-recovery), tariffs are adjusted upward to compensate.

\textit{Capacity allocation} is the competitive bidding process for access to capacity rights, which has historically transitioned from a first-come, first-served (FCFS) system to a regulated competitive bidding process under the EU's Third Energy Package \citep{EuropeanCommission2009ThirdPackage}. The modern mechanism involves auctioning capacity products by duration: annual capacity is auctioned once per year; quarterly capacity is auctioned 1-4 times per quarter; monthly capacity is auctioned once per month; daily capacity is auctioned once per day; and within-day auctions occur multiple times per day.
Auction formats vary by product type: daily and within-day auctions use uniform pricing, in which all successful bids pay the same clearing price, whereas annual, quarterly, and monthly capacity products typically use an ascending-clock auction format. In the ascending clock auction, prices start at the base rate and increase incrementally over multiple rounds. This format allows participants to adjust their bids incrementally until the supply-and-demand equilibrium is reached.

\subsection{Nomination, Renomination, Allocation, and Physical Flow}
Shippers submit their \textit{nomination}, which is the schedule of hourly volumes (in MWh/h) they intend to inject or withdraw at designated entry and exit points after booking capacity. The European gas system operates on a standardised 06:00 to 06:00 CET as displayed in Figure \ref{fig:gas_market_timeline}, due to historical offtake patterns typically starting early in the morning. Nominations for gas day D must be submitted by 14:00 CET on D-1 and validated by TSOs between 14:00 and 18:00. Day-ahead and within-day capacity is booked after the first nomination rounds. 

\textit{Renomination} allows shippers to adjust their renominations to adapt to market conditions or system constraints. This process opens after 18:00 on D-1, allowing hourly adjustments with a two-hour lead time. For example, before 10:00 on gas day D-1, a renomination can be placed for all hours from 12:00 to 06:00, covering the remaining hours of the gas day. Renominations can cover the entire 24-hour gas day until 04:00 on gas day D. After 04:00, the renomination window shortens by one hour for each ensuing hour, with the final renomination window at 03:00 on gas day D, covering from 05:00 to 06:00.

\begin{figure}
\centering
\begin{tikzpicture}[
    scale=1,
    timeline/.style={->, >=stealth, very thick},
    tick/.style={thick},
    time label/.style={font=\small\bfseries},
    phase box/.style={draw, dashed, rounded corners=3pt, fill=gray!10, 
                     minimum height=0.8cm, align=center, font=\footnotesize},
    day label/.style={font=\small\bfseries, red!70!black}
]

% Timeline arrow
\draw[timeline] (0,0) -- (15,0);

% Day markers and labels

\draw[decorate, decoration={brace, amplitude=6pt, mirror}] 
    (0, -0.3) -- (10, -.3);
    
\node[box, draw=blue!60!black, rounded corners, fill=blue!10, minimum width=2.5cm] at (5, -1) {
    \textbf{Gas day ahead D-1}
};

\draw[decorate, decoration={brace, amplitude=6pt, mirror}] 
    (10, -0.3) -- (15, -.3);
    
\node[box, draw=green!60!black, fill=green!10, rounded corners, minimum width=2.5cm] at (12.5, -1) {
    \textbf{Gas day D}
};

% Time ticks and labels
\foreach \x/\time in {0/6:00, 1.6/10:00, 3.3/14:00, 5/18:00, 6.6/22:00, 8.3/2:00, 10/6:00, 15/6:00} {
    \draw[tick] (\x,-0.1) -- (\x,0.1);
    \node[time label] at (\x,0.4) {\time};
}

% Phase boxes with improved positioning
%\node[phase box, draw=blue!60!black, fill=blue!10, minimum width=2.5cm] at (3.3, 1.5) {
%     \shortstack{\textbf{Nomination}\\\textbf{Review}}
%};

\node[phase box, draw=blue!60!black, fill=blue!10, minimum width=2.5cm] at (3.3, 1.5) {
    \shortstack{Nominations\\submission}
};

\node[phase box, draw=blue!60!black, fill=blue!10, minimum width=2.5cm] at (6,1.5) {
    \shortstack{Renominations\\start}
};

\node[phase box, draw=green!60!black, fill=green!10, minimum width=2.5cm] at (10, 1.5) {
    Gas day starts
};

\node[phase box, draw=green!60!black, fill=green!10, minimum width=2.5cm] at (15, 1.5) {
    Gas day end
};

% Connecting lines from boxes to timeline
\draw[gray, dashed] (3.3, 1) -- (3.3, 0.6);
\draw[gray, dashed] (5, 0.6) -- (5, 0.8) -- (6, 0.8) -- (6, 1);

%Gas day start
\draw[gray, dashed] (10, 1.1) -- (10, 0.6);
% Gas day end
\draw[gray, dashed] (15, 1.1) -- (15, 0.6);

% Timeline periods with colored backgrounds
\fill[blue!20, opacity=0.3] (0,-0.2) rectangle (10,0.2);
\fill[green!20, opacity=0.3] (10,-0.2) rectangle (15,0.2);

% Legend
%\node[anchor=north west, font=\footnotesize] at (0,3.2) {
%    \begin{tabular}{@{}l@{\hspace{0.5em}}l@{}}
%    \textcolor{green!70!black}{\rule{0.3cm}{0.1cm}} & Nomination Submission \\
%    \textcolor{orange!70!black}{\rule{0.3cm}{0.1cm}} & Nomination Review \\
%    \textcolor{blue!70!black}{\rule{0.3cm}{0.1cm}} & Renomination Period \\
%    \end{tabular}
%};

\end{tikzpicture}
\caption{\textbf{EU Gas Day-Ahead and Day-Ahead Gate Closures.} Schematic representation of the gas day-ahead market timeline showing key operational phases}
\label{fig:gas_market_timeline}
\end{figure}

\textit{Allocation} is performed by the TSO at the end of each gas day. This accounting process assigns the actual metered gas volumes to individual shippers based on their nominations and the physical flow. If the physical flow matches the nomination, allocation equals nomination. However, if deviations occur due to renomination, network constraints, or metering errors, the renomination location reflects the reconciled volume actually delivered or received. Generally, at border points, allocation usually matches the nomination, while at consumption points, it is based on metering. Less frequent cases of modification to allocation accounting include the use of interruptible capacity or reductions in allocations due to unexpected maintenance.

\textit{Physical flow} refers to the actual measured amount of gas transported and consumed, and it does not have to match the nomination or allocation exactly in real time. The system is optimised based on factors such as coordination with counterparts and operational constraints. Once nominations and renominations are in place, gas flows according to the updated schedules. 

In Figure \ref{fig:capacity-allocation}, the lower plot illustrates these dynamics. While the nominations are steady throughout the displayed period, renominations differ from the original schedule. These adjustments are reflected in the blue line (allocation), which updates accordingly. The green line shows the physical flow, which may deviate slightly but generally follows the contractual path.

\textit{Portfolios} are associated with each shipper and serve as formal documentation of their booked capacities at various points, nomination schedules, and resulting allocations across all points where they hold capacity rights. In an entry-exit system, TSOs monitor the net position of each shipper's portfolio, which is the sum of all nominations across the portfolio. For example, if a shipper nominates +10 MWh at entry point A and -10 MWh at exit point B, the positions at each point are +10 MWh and -10 MWh, respectively, while the net position of the entire portfolio remains 0 MWh. Shippers must ensure that the portfolio's net position is balanced, meaning that nominations at entry points must correspond to matching exit points on the same day (conditional on WDOs). If a supply constraint is imposed on the shipper, requiring a change from 10 MWh to 9 MWh at entry point A, the shipper must renominate the exit volume at exit point B to 9 MWh to maintain a balanced portfolio.

\subsection{Balancing framework}
\label{subsec:balancing}
\begin{figure}
\centering
\begin{tikzpicture}[scale=1]
% Nature journal color scheme
\definecolor{natureblue}{RGB}{0, 119, 181}      % Nature primary blue
\definecolor{natureteal}{RGB}{0, 150, 136}      % Nature teal
\definecolor{natureorange}{RGB}{255, 87, 34}    % Nature orange
\definecolor{naturegreen}{RGB}{76, 175, 80}     % Nature green
\definecolor{naturegray}{RGB}{97, 97, 97}       % Nature dark gray
\definecolor{axiscolor}{RGB}{33, 33, 33}        % Deep gray for axes
\definecolor{neutralbrown}{RGB}{121, 85, 72}    % Neutral brown for support linepack

% --- Linepack Bar (Left side) ---
% Left brace for total linepack
\draw[thick, decorate, decoration={brace, amplitude=6pt, mirror}] 
    (-1.25, 1.5) -- (-1.25, -3);
\node[anchor=east, font=\footnotesize\bfseries] at (-1.5, -0.75) {Total Linepack};

% Pressure level indicators
\draw[thick] (-0.3, 1.5) -- (-0.2, 1.5);
\draw[thick] (-0.3, -2.0) -- (-0.2, -2.0);
\node[anchor=west, font=\footnotesize] at (-1.2, 1.5) {$P_{\text{max}}$};
\node[anchor=west, font=\footnotesize] at (-1.2, -2.0) {$P_{\text{min}}$};

% --- Main Linepack Regions ---
% Flow Linepack (triangular region)
\fill[natureblue, opacity=0.7] (0, 1.5) -- (6, 0.5) -- (0, 0.5) -- cycle;
\draw[thick, color=natureblue] (0, 1.5) -- (6, 0.5) -- (0, 0.5) -- cycle;
\node[font=\footnotesize\bfseries, color=white] at (1.75, .8) {Flow Linepack};

% System Linepack Flexibility (rectangular region)
\fill[natureteal, opacity=0.7] (0, 0.5) rectangle (6, -1.5);
\draw[thick, color=natureteal] (0, 0.5) rectangle (6, -1.5);
\node[font=\footnotesize\bfseries, color=white] at (3.0, -0.5) {System Linepack Flexibility};

% Reserve Linepack (rectangular region)
\fill[natureorange, opacity=0.7] (0, -1.5) rectangle (6, -2.0);
\draw[thick, color=natureorange] (0, -1.5) rectangle (6, -2.0);
\node[font=\footnotesize\bfseries, color=white] at (3.0, -1.75) {Reserve Linepack};

% Support Linepack (rectangular region)
\fill[neutralbrown, opacity=0.7] (0, -2.0) rectangle (6, -3.0);
\draw[thick, color=neutralbrown] (0, -2.0) rectangle (6, -3.0);
\node[font=\footnotesize\bfseries, color=white] at (3.0, -2.5) {Support Linepack};

% --- Coordinate System ---
% Axes with professional styling
\draw[->, thick, color=axiscolor] (0, -3.2) -- (6.8, -3.2) 
    node[anchor=north west, font=\footnotesize] {Pipeline Length [km]};
\draw[->, thick, color=axiscolor] (-0.2, -3) -- (-0.2, 1.8) 
    node[anchor=south east, font=\footnotesize] {Pressure [bar]};

% --- Optional: Add dimension arrows for clarity ---
\draw[<->, thick, color=naturegray] (6.2, 1.5) -- (6.2, 0.5);
\node[anchor=west, font=\scriptsize, color=naturegray] at (6.3, 1.0) {$\Delta P^{\text{Flow}}$};

\draw[<->, thick, color=naturegray] (6.2, 0.5) -- (6.2, -1.5);
\node[anchor=west, font=\scriptsize, color=naturegray] at (6.3, -0.5) {$\Delta P^{\text{System}}$};

\draw[<->, thick, color=naturegray] (6.2, -1.5) -- (6.2, -2.0);
\node[anchor=west, font=\scriptsize, color=naturegray] at (6.3, -1.75) {$\Delta P^{\text{Reserve}}$};

\draw[<->, thick, color=naturegray] (6.2, -2.0) -- (6.2, -3.0);
\node[anchor=west, font=\scriptsize, color=naturegray] at (6.3, -2.5) {$\Delta P^{\text{Support}}$};

\end{tikzpicture}
\caption{\textbf{Linepack flexibility regions in natural gas pipeline systems.} The diagram illustrates the pressure-dependent allocation of linepack capacity across different operational zones, from maximum system pressure ($P_{\text{max}}$) to minimum operational pressure ($P_{\text{min}}$). Illustration based on \cite{Energinet2024ProposalNetwork}.}
\label{fig:linepack_flexibility}
\end{figure}
In an entry-exit system, the TSO manages the balance of its balancing zone by monitoring the net position of each shipper's portfolio. A portfolio's net position is long when it injects more gas than it withdraws (surplus) and short when it withdraws more gas than it injects (deficit) within a given timeframe. The \textit{System Balancing Signal} (SBS) represents the aggregate net position of all portfolios within a balancing zone, reflects the overall system balance, and is therefore expressed in MWh. The SBS can be in three states: positive SBS indicates a long system with excess gas and higher-than-nominal pressure; negative SBS indicates a short system with a gas deficit and lower pressure; and zero SBS indicates a balanced system with no net imbalance.

Portfolios are classified based on their contribution to the system imbalance. Shippers whose net positions align with the direction of the imbalance are termed \textit{causers}, while those whose net positions counteract the imbalance are termed \textit{helpers}. When the system is long, long portfolios are causers and short portfolios are helpers, and conversely, when the system is short, short portfolios are causers and long portfolios are helpers.

From a physical perspective, a gas system is imbalanced when the pipeline pressure falls outside safety limits. To link the financial balancing of the portfolios to the physical pipeline balancing, the TSO introduced the concept of \textit{linepack}. Linepack corresponds to the pressure of the gas stored in pipelines expressed in MWh and, according to Figure \ref{fig:linepack_flexibility}, is categorised into components: support linepack ensures minimum pressure levels, reserve linepack is allocated for TSO balancing actions and system survival, system linepack flexibility is available for flexible use by network users, and flow linepack enables hydrogen transport through pressure differences \citep{Energinet2024ProposalNetwork}.

The principles of the balancing scheme are outlined in the European Balancing Network Code, which regulates how TSOs should establish a balancing regime for their zones \citep{Entsog2018BalancingOverview}. The end goal of a functioning balancing regime is to allow shippers to collectively balance their portfolios and minimise TSO costs through \textit{balancing actions}, defined as interventions that occur when the SBS moves outside predefined safety boundaries and require the TSO to buy or sell gas on the spot market. For systems exceeding limits, TSOs sell excess gas; for systems falling short, they buy gas. Since balancing actions are conducted at marginal spot market prices, they can be costly, with the costs distributed proportionally across the causer portfolios.

A common principle adopted among EU TSOs for evaluating the need for balancing action is the use of \textit{buffer zones}, which divide the system linepack flexibility into zones. This is applied differently among EU countries; for example, the Dutch TSO Gasunie organises the following \citep{GasunieTransportServices2025SBSPOS}: the green zone operates with no balancing action as long as the SBS remains within limits; the yellow zone triggers end-of-day balancing actions under certain circumstances; the orange zone requires next-hour balancing actions when the SBS enters or remains in this zone; and the red zone demands immediate next-hour balancing actions. 

\begin{figure}
\centering
\includegraphics[width=\linewidth]{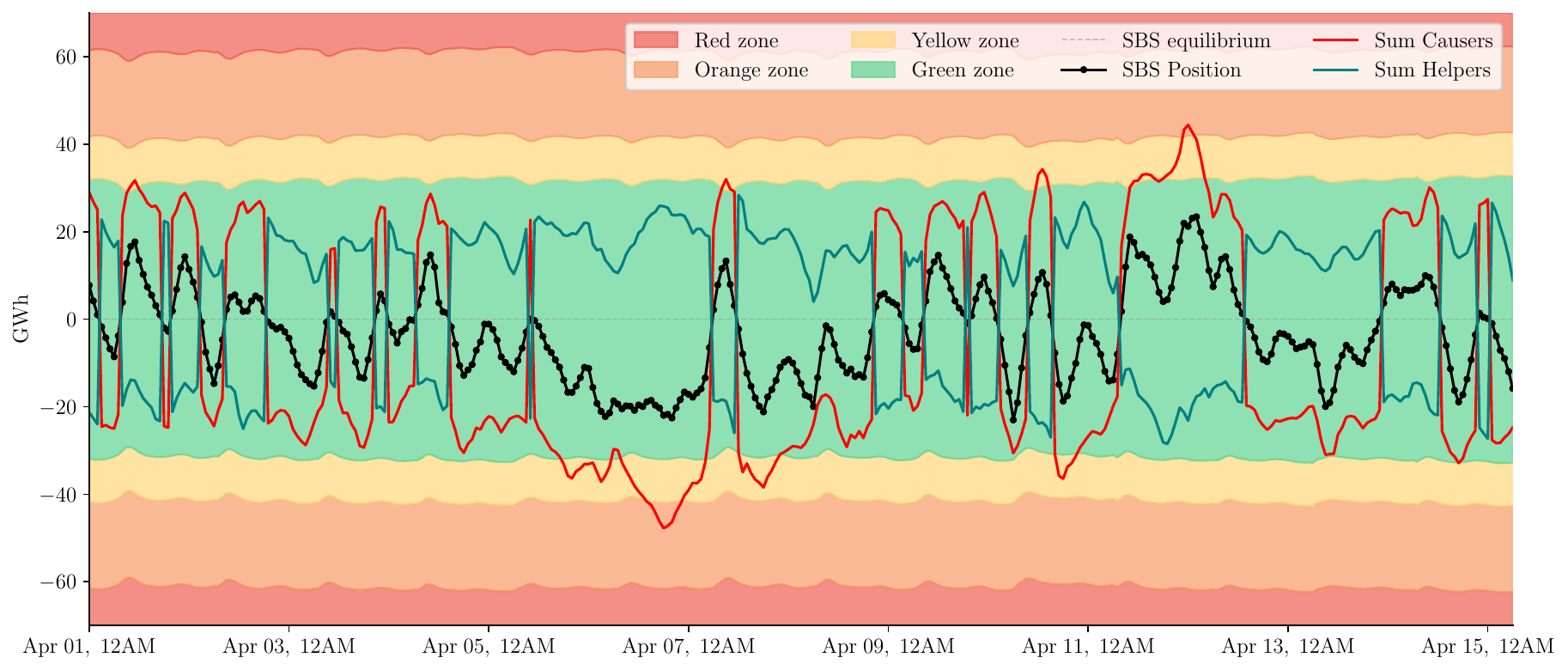}
\caption{\textbf{Example of System Balancing Signal (SBS) and portfolio positions}. Data displayed for the Dutch balancing zone, from April $1^\text{th}$ to $15^\text{th}$ 2025 (\cite{GasunieTransportServices2025SBSPOS}).}
\label{fig:SBS_sample}
\end{figure}

Figure \ref{fig:SBS_sample} illustrates the SBS evolution in the Dutch network over two weeks. On April $7^\text{th}$ 2025, the SBS briefly changed from negative to positive and back again, temporarily reversing the roles of causers and helpers as the system fluctuated between long and short conditions. No balancing action was needed during this period because the SBS remained within the acceptable (green) zone.

\textit{Virtual Trading Points} (VTPs) are non-physical entry-exit points within a balancing zone that allow gas to be traded between different market participants, such as producers and consumers, without requiring a physical gas exchange. VTPs are essential for enabling gas ownership transfers; without them, there is no basis for a market to function. Consider the example presented in Figure \ref{fig:portfolio_balancing_table} of a balancing zone with shippers A and B: if Shipper A wants to sell gas to Shipper B, Shipper A injects 10 kWh of gas at an entry point, and Shipper B withdraws 10 kWh at an exit point, completing the physical gas exchange. Although the zone's SBS is 0 kWh, individual portfolios are imbalanced, as shown in the left plot of the figure: Shipper A's portfolio is long by 10 kWh, and Shipper B's portfolio is short by 10 kWh. To balance these portfolios without a VTP, Shipper A would need to purchase capacity at an exit point, physically exit and transfer the gas to Shipper B, who would then need to buy capacity at an entry point before reinjecting and withdrawing the gas. This process is operationally and financially inefficient due to additional capacity costs and the physical handling of the gas. The purpose of VTPs is to enable different entities (e.g., producers and consumers) to exchange gas. 

In the right plot of Figure \ref{fig:portfolio_balancing_table}, the VTP allows shippers to transfer gas ownership virtually within the system, avoiding physical gas movement. As in the previous instance, Shipper A nominates the injection of gas, and Shipper B nominates the withdrawal. Additionally, both place a nomination at the VTP: Shipper A virtually exits the gas, and Shipper B virtually enters the gas, without any physical exchange. This virtual nomination balances the portfolios, ensuring the gas enters and exits the system only once. In reality, shippers submitting nominations to the VTP must specify trade amounts, delivery periods, and trade directions (entry or exit). For OTC trades, the TSO matches opposing nominations to complete the virtual transaction, whereas in exchange-based trades, matching is anonymous through the trading platform.

\begin{figure}
\centering
\begin{tikzpicture}[
    scale=1,
    % Box styles
    box/.style={draw, rounded corners=3pt, minimum width=2.2cm, minimum height=0.8cm, 
                align=center, font=\footnotesize},
    header box/.style={box, draw=blue!60!black, fill=blue!10, font=\footnotesize\bfseries},
    blue box/.style={box, draw=blue!60!black, fill=blue!5},
    green box/.style={box, draw=green!60!black, fill=green!5},
    red box/.style={box, draw=red!60!black, fill=red!5},
    teal box/.style={box, draw=teal!60!black, fill=teal!5},
    gray box/.style={box, draw=gray!60!black, fill=gray!10, minimum width=2.5cm, font=\footnotesize\bfseries},
    % Text styles
    positive/.style={color=green!70!black, font=\footnotesize\bfseries},
    negative/.style={color=red!70!black, font=\footnotesize\bfseries},
    neutral/.style={color=gray!70!black, font=\footnotesize\bfseries},
    % Arrow styles
    trade arrow/.style={->, >=stealth, thick, green!70!black},
    % Section titles
    section title/.style={font=\large\bfseries, blue!80!black}
]

% Section A: Without Virtual Trading Point

% Table headers for section A
\node[blue box] at (1.5,6) {\textbf{Shipper A}};
\node[teal box] at (4 ,6) {\textbf{Shipper B}};

% Row 1: Nomination
\node[gray box] at (-1,5) {Nomination};
\node[blue box] at (1.5,5) {-10 kWh};
\node[teal box] at (4,5) {+10 kWh};

% Row 2: VTP trade (empty for section A)
\node[gray box] at (-1,4) {VTP trade};
\node[blue box] at (1.5,4) {---};
\node[teal box] at (4,4) {---};

\draw[gray, dashed] (-2.2, 3.5) -- (5.1, 3.5);
% Row 3: Net position
\node[gray box] at (-1,3) {Net position};
\node[red box] at (1.5,3) {-10 kWh};
\node[red box] at (4,3) {+10 kWh};

% Section B: With VTP

% Table headers for section B
\node[blue box] at (9.5,6) {\textbf{Shipper A}};
\node[teal box] at (12,6) {\textbf{Shipper B}};

% Row 1: Nomination
\node[gray box] at (7,5) {Nomination};
\node[blue box] at (9.5,5) {-10 kWh};
\node[teal box] at (12,5) {+10 kWh};

% Row 2: VTP trade
\node[gray box] at (7,4) {VTP trade};
\node[blue box] at (9.5,4) {+10 kWh\\{\tiny (VTP)}};
\node[teal box] at (12,4) {-10 kWh\\{\tiny (VTP)}};
\draw[gray, dashed] (5.9, 3.5) -- (13.1, 3.5);

% Row 3: Net position
\node[gray box] at (7,3) {Net position};
\node[green box] at (9.5,3) {0 kWh};
\node[green box] at (12,3) {0 kWh};

% Trade arrows
\draw[trade arrow] (10.4,4.2) -- (11.2,4.2);
\draw[trade arrow] (11.2,3.8) -- (10.4,3.8);

% Summary result boxes
\node[draw=red!60!black, rounded corners=3pt, fill=red!15, font=\footnotesize, 
      align=center, minimum width=2.5cm] at (2.75, 2) {
    Portfolios \textcolor{red!70!black}{\textbf{imbalanced}}
};

\node[draw=green!60!black, rounded corners=3pt, fill=green!15, font=\footnotesize, 
      align=center, minimum width=2.5cm] at (10.75, 2) {
    Portfolios \textcolor{green!70!black}{\textbf{balanced}}
};

\draw[gray, opacity=0.5, line width=0.5pt] (5.4, 1.5) -- (5.4, 6.8);

\end{tikzpicture}
\caption{\textbf{Comparison of portfolio balancing mechanisms without (left plot) and with (right plot) Virtual Trading Points (VTP)}. Without VTP, physical gas movement is required to balance portfolios; with VTP, gas ownership transfers are virtual, eliminating the need for physical movement and improving operational efficiency.}
\label{fig:portfolio_balancing_table}
\end{figure}

\subsection{Commodity trading and gas spot market}
Natural gas in Europe is traded through two main mechanisms: over-the-counter (OTC) bilateral contracts and centralised forms, both of which rely on virtual trading points to transfer gas ownership without physical movement. Before gas market liberalisation, companies were vertically integrated, supplying the commodity to the customers and managing transportation infrastructure. The EU Third Energy Package \citep{EuropeanCommission2009ThirdPackage} enforced unbundling, making TSOs responsible for infrastructure and allowing other parties to trade gas as a commodity.

\textit{OTC Contracts} facilitate bilateral trades between buyers and sellers and are typically characterised by long-term duration and customisable pcustomisable schemes. Unlike standardised contracts, OTC contracts are fully customised by parties, allowing agreement on volume, duration, delivery profiles, and pricing structures. This flexibility is essential for long-term agreements (e.g., 30-year contracts) that standardised contracts cannot accommodate. Historically, OTC contracts have dominated gas markets because suppliers needed long-term agreements to secure bankability. However, because OTC contracts are bilateral and not publicly traded, they lack transparency in pricing and trading volumes, with external agencies reporting prices through broker and market participant surveys.

\textit{Exchange trading} involves standardised platforms such as ICE (Intercontinental Exchange) and EEX (European Energy Exchange), which operate across both the US and European markets. These exchanges offer standardised contracts, including futures with delivery periods ranging from monthly to yearly. Unlike the electricity spot market, which clears hourly at a single price, the gas spot market operates as a continuous trading system with centralised participants who can execute transactions at any time during trading hours, resulting in a continuous set of prices throughout the day. Exchanges calculate volume-weighted average prices at market close to determine the daily settlement price, commonly referred to as the "spot price." Initially, OTC contracts had the highest liquidity, but over time, exchanges have increased their liquidity.

\textit{The role of VTP} is essential in ensuring the correct transfer of gas ownership between portfolios. Most gas trades, whether on exchanges or in OTC markets, use the VTP to transfer gas ownership among portfolios. For example, if a trade occurs between two shippers on EEX or ICE, the trade information is sent to the TSO, which then enters it into the shippers' portfolios via a nomination at the VTP. Exchange trades integrate automatically into nomination systems, while OTC transactions require manual or broker-supported nomination processes. VTP liquidity is the volume of trades executed at the VTP for that balancing zone. The largest VTP in Europe by liquidity is the Dutch Title Transfer Facility (TTF), followed by the British National Balancing Point (NBP). Although the term "VTP price" is commonly used to refer to the natural gas price, it is technically incorrect. The correct terminology is "EEX/ICE price of gas traded, settled at the TTF," although "TTF price" is widely accepted. The VTP price reflects volume-weighted averages of daily trades for specific products. As the most liquid VTP, the TTF serves as an effective continental price proxy, with its price as the reference, as seen in Figure \ref{fig:TTF-prices}. However, network congestion across balancing zones can lead to significant regional price deviations.

\begin{figure}
\centering
\includegraphics[width=\linewidth]{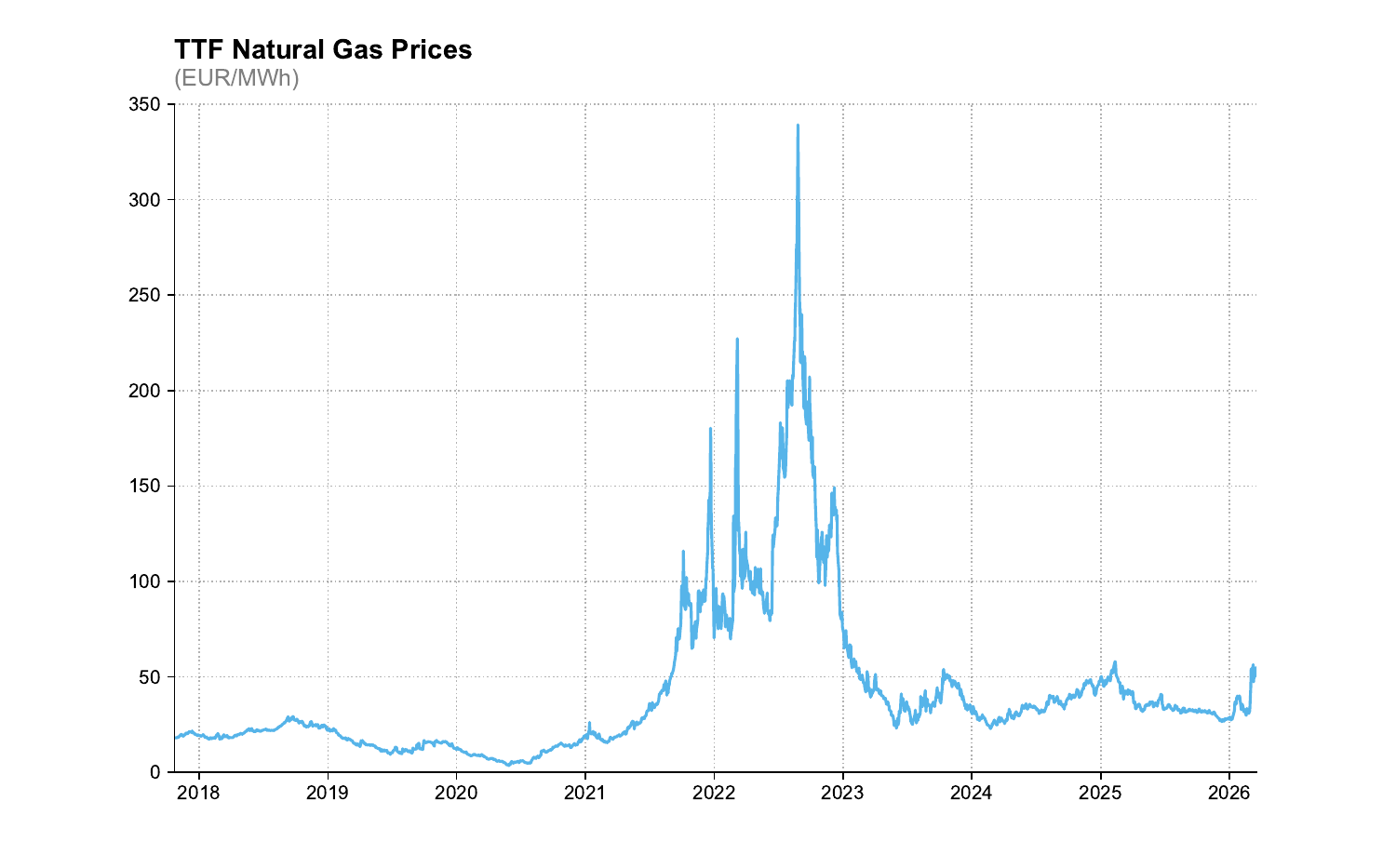}
\caption{\textbf{Historical gas prices in the Dutch gas balancing zone}. Prices are referred to as the TTF, a proxy for VTP in Europe \citep{YahooFinance2025DutchData}.}
\label{fig:TTF-prices}
\end{figure}

\section{EU Energy Packages}
\label{sec:energy_packages}
The EU transformed its energy markets from monopolistic structures in the 1990s to competitive markets through five progressive Energy Packages:

\begin{enumerate}
    \item \textbf{First Energy Package (1996-1998)}
Initiated energy market liberalisation with third-party network access and unbundling of Transmission System Operators (TSOs), separating competitive activities (production/supply) from natural monopolies (transmission/distribution).

Legislation: Electricity Directive 96/92/EC (1996), Gas Directive 98/30/EC (1998)

\item \textbf{Second Energy Package (2003)}
Advanced liberalisation by enabling consumer choice of energy suppliers, requiring legal unbundling of TSOs, and establishing independent National Regulatory Authorities (NRAs).

Legislation: Electricity Directive 2003/54/EC, Gas Directive 2003/55/EC, Regulation 1228/2003 (cross-border electricity access)

\item \textbf{Third Energy Package (2009)}
Further enhanced TSO unbundling requirements, moving beyond legal separation to ensure complete independence of supply and generation activities from network operations. This addressed the limitation that early market liberalisation allowed trading only at network entry and exit points, hindering efficient commodity exchange.

The package strengthened NRA independence, established ACER, and created ENTSO-E and ENTSO-G for cross-border coordination. It introduced EU network codes, including the Gas Balancing NC, the Capacity Allocation Mechanisms NC, the Congestion Management Procedures, and the Harmonised Transmission Tariff Structures NC.

\textbf{Key Innovation}: VTPs, as an administrative concept, enable gas commodity exchanges within the network system. These virtual points allow parties to trade gas while it remains in the system, facilitating the development of competitive exchanges and OTC contracts.

Legislation: Electricity Directive 2009/72/EC, Gas Directive 2009/73/EC, ACER Regulation 713/2009/EC, Electricity Regulation 714/2009/EC, Gas Regulation 715/2009/EC

\item \textbf{Fourth Energy Package / "Clean Energy Package" (2019)}
Modernised the EU energy system with consumer empowerment tools and established 2030 climate targets: 40\% GHG emission reduction, 32\% renewables, 32.5\% energy efficiency.

Legislation: Four directives (Electricity 2019/944, RED II 2018/2001, Energy Efficiency 2018/2002, Buildings Performance 2018/844) and four regulations (Electricity 2019/943, Risk Preparedness 2019/941, ACER 2019/942, Governance 2018/1999)

\item \textbf{Fifth Energy Package (2024)}
Aligned energy policy with climate neutrality goals and extended the gas framework to hydrogen. Incorporated the REPowerEU plan to accelerate the Russian fossil fuel phase-out, enhance energy diversification, and reform electricity markets with stronger long-term contracts and consumer protections.

\end{enumerate}

\section{RFNBO Legislation}
\label{sec:rfnbo_legislation}

The Renewable Energy Directive (RED) establishes the EU's framework for the development of clean energy across all sectors, aiming to achieve climate neutrality by 2050. As part of its decarbonisation strategy, the EU aims to develop Renewable Fuels of Non-Biological Origin (RFNBO), such as renewable hydrogen, which is essential for decarbonising energy-intensive industrial processes and the transport sector where direct electrification is not feasible.

To produce RFNBO, strict criteria must be met to ensure that the electricity used is renewable. The Delegated Regulation (EU) 2023/1184 outlines detailed rules for determining when electricity can be considered fully renewable:

\begin{itemize}
    \item \textbf{Electricity Sourced via Direct Connection:}
    \begin{enumerate}
        \item Electricity obtained directly from a renewable energy installation connected to the RFNBO production site qualifies as fully renewable without additional requirements.
    \end{enumerate}
    
    \item \textbf{Electricity Sourced from the Grid:}
    When the electricity for RFNBO production is sourced from the grid, several stringent conditions must be met:
    \begin{enumerate}
        \item \textbf{Additionality:} The renewable energy installation must have commenced operation not earlier than 36 months before the electrolyser and must not have received public financial support. This ensures that the renewable electricity is additional and not diverting existing resources.
        \item \textbf{Geographical Correlation:} The renewable electricity source and the hydrogen production site must be located in the same bidding zone to avoid grid congestion and ensure efficient energy transfer.
        \item \textbf{Temporal Correlation:} Until 2029, the electricity consumption and production must correlate within the same month to ensure the renewable electricity is actively contributing to hydrogen production. From 2030, this correlation must occur hourly, aligning more closely with real-time production and consumption.
    \end{enumerate}
\end{itemize}

The Delegated Regulation (EU) 2023/1185 complements these criteria by establishing a methodology for calculating the lifecycle greenhouse gas emissions associated with RFNBOs, ensuring comprehensive sustainability.

Despite these regulations, there are criticisms of the strictness of the criteria, which may limit the rapid deployment of renewable hydrogen production needed to achieve the EU's climate goals. These concerns highlight the need for balancing stringent sustainability requirements with practical implementation to foster the growth of a viable hydrogen market.
\section{Gas Congestion Management}

Gas congestion occurs in two forms: commercial congestion, when capacity is booked but not utilised, and physical congestion, when infrastructure constraints prevent the required gas flows. Commercial congestion can be addressed through reallocation mechanisms, while physical congestion requires infrastructure expansion.

Congestion typically manifests when significant price differences exist between gas hubs, creating arbitrage opportunities that drive shippers to rapidly book transport capacity. This phenomenon is particularly pronounced in short-term markets, such as day-ahead or within-day auctions. Physical constraints may exacerbate congestion during extreme weather events when supply limitations create security concerns.TSOs employ two categories of congestion management tools: direct or indirect.

\textbf{Direct tools} are specifically designed for congestion management:
\begin{enumerate}
    \item \textbf{Use It or Lose It (UIOLI)} enables regulators to reallocate unused booked capacity. However, UIOLI implementation remains rare in Europe due to complex regulatory requirements involving multiple National Regulatory Authorities and legal concerns about the reallocation of capacity from private entities.
\end{enumerate}

\textbf{Indirect tools} serve other primary purposes while alleviating congestion as a secondary effect:
\begin{enumerate}
    \item \textbf{Auctions} create pricing mechanisms that impose premiums during congestion periods, making inefficient capacity use costly and ensuring fair allocation.
    
    \item \textbf{Counterflow} optimises capacity utilisation by enabling commercial transactions in both directions at bidirectional connection points. While gas flows physically in only one direction at any given time, commercial capacity can be sold bidirectionally. This allows offsetting flows (e.g., 100 units north and 100 units south resulting in zero net physical flow) while maintaining commercial capacity utilisation.
    
    \item \textbf{Overnomination} provides short-term interruptible capacity that can be nominated up to two hours before the gas hour. Payment occurs only upon nomination, and this capacity receives the lowest priority in the allocation hierarchy: (1) Firm capacity, (2) Interruptible capacity, (3) Within-day interruptible capacity, (4) Overnomination. This mechanism maximises available capacity utilisation while providing flexibility for last-minute demand adjustments.
\end{enumerate}

The combination of these tools, particularly auction-based allocation with appropriate pricing mechanisms, has proven effective in reducing congestion while maintaining efficient capacity utilisation across European gas networks.

\end{document}